\documentstyle[amssymb,12pt]{article}
\title{Borel completeness of some $\aleph_0$-stable theories}
\author{Michael C.\
Laskowski\thanks{Partially supported
by NSF grants DMS-0600217 and 0901334.}\\
Department of Mathematics\\University of Maryland
\and
S. Shelah\thanks{
Partially supported by U.S.-Israel Binational Science Foundation
Grant no.\ 2002323
and Israel Science Foundation Grant no.\ 242/03.
 Publication no.\ 1016.
Both authors were partially supported by NSF grants DMS-0600940 and 1101597.
}
\\
Department of Mathematics\\
Hebrew University of Jerusalem\\
Department of Mathematics\\
Rutgers University
}

\newbox\smilebox
\newbox\anchorbox
\newbox\noanchorbox
\newbox\tempbox

\setbox\smilebox=\hbox{$\smile$}

\def\anchor{\hbox{\vtop{
           \hbox to \wd\smilebox{\hfil\vrule width.4pt height7pt depth1pt\hfil}
           \vskip  -11.5truept
           \hbox to \wd\smilebox{\hfil$\smile$\hfil}}}}
\setbox\anchorbox=\anchor
\def\noanchor{\hbox{\vtop{
           \hbox to \wd\anchorbox{\hfil\anchor\hfil}
           \vskip -14truept
           \hbox to \wd\anchorbox{\hfil/\hfil}}}}
\setbox\noanchorbox=\noanchor

\def\fg#1#2#3{\setbox\tempbox=\hbox{$\scriptstyle{#2}$}
\ifnum\wd\anchorbox>\wd\tempbox\dimen255=\wd\anchorbox
\else\dimen255=\wd\tempbox\fi
{#1\,\vtop{\hbox to \dimen255{\hfil\anchor\hfil}
           \vskip -6truept
           \hbox to \dimen255{\hfil$\scriptstyle{#2}$\hfil}}
           \,#3}}

\def\nfg#1#2#3{\setbox\tempbox=\hbox{$\scriptstyle{#2}$}
\ifnum\wd\noanchorbox>\wd\tempbox\dimen255=\wd\noanchorbox
\else\dimen255=\wd\tempbox\fi
{#1\,\vtop{\hbox to \dimen255{\hfil\noanchor\hfil}
           \vskip -6truept
           \hbox to \dimen255{\hfil$\scriptstyle{#2}$\hfil}}
           \,#3}}

\setbox1=\hbox{$\bot$}

\def\north#1#2{#1\,
\hbox{$\bot$\llap {\hbox to\wd1 {\hfil $/$\hfil}}}
\,#2}

\def\nao#1#2#3{#1\  \hbox{\vtop{ 
\baselineskip=4pt
\hbox{$\bot$\llap {\hbox to\wd1 {\hfil $/$\hfil}}
\hskip .05em \llap{\hbox{$^{\scriptscriptstyle{a}}$}}}\hbox{$\scriptstyle
{#2}$}}}\, #3}

\def\bp{\par{\bf Proof.}$\ \ $}

\def\includeE#1{{\lhook\kern-3.5pt\joinrel\smash{
    \mathop{\longrightarrow}\limits^{#1}}}}

\def\efor/{Example~\ref{E4}}

\def\BL/{Baldwin--Lachlan}
\def\Bu/{Buechler}
\def\Hr/{Hrushovski}
\def\lm/{locally modular}
\def\wm/{weakly minimal}
\def\nm/{non--modular}
\def\ss/{superstable}
\def\ud/{unidimensional}
\def\sm/{strongly minimal}

\def\abar{\overline{a}}

\def\hbar{\overline{h}}

\def\xbar{\overline{x}}

\def\Mbar{\overline{M}}

\def\dom{{\rm dom}}

\def\tp{{\rm tp}}

\def\tr/{trivial}
\def\nt/{non--trivial}
\def\st/{strong type}

\def\TV/{Tarski--Vaught}

\def\sc/{sound construction}
\def\ac/{atomic construction}

\def\fal/{functional}
\def\upl/{unique parallel lines}
\def\chp/{categorical in a higher power}

\def\conc{{\char'136}}
\def\abar{\bar{a}}

\def\phi{\varphi}
\def\C{{\frak  C}}

\def\B{{\cal B}}
\def\C{{\frak  C}}
\def\CC{{\cal C}}
\def\U{{\cal U}}
\def\diff{{\rm diff}}
\def\D{{\cal D}}
\def\E{{\cal E}}
\def\F{{\cal F}}
\def\Lplusminus{{L^{\pm}}}
\def\M{{\cal M}}
\def\P{{\mathbf P}}
\def\Pcal{{\cal P}}
\def\T{{\cal T}}

\def\d{{\mathfrak d}}

\def\R{{\cal R}}

\def\S{{\cal S}}
\def\F{{\cal F}}
\def\T{{\cal T}}

\def\alphabar{\overline{\alpha}}
\def\betabar{\overline{\beta}}
\def\bfequiv{\equiv_{\infty,\aleph_0}}
\def\range{{\rm range}}
\def\tp{{\rm tp}}

\def\wt{{\rm wt}}

\def\dom{{\rm dom}}

\def\bp{{\bf Proof.}\quad}
\def\endproof{\medskip}
\def\<{\langle}
\def\>{\rangle}

\newtheorem{Theorem}{Theorem}[section]
\newtheorem{Proposition}[Theorem]{Proposition}
\newtheorem{Definition}[Theorem]{Definition}

\newtheorem{Lemma}[Theorem]{Lemma}
\newtheorem{Corollary}[Theorem]{Corollary}
\newtheorem{Fact}[Theorem]{Fact}
\begin{document}
\maketitle
\begin{abstract}
We study $\aleph_0$-stable theories, and prove that
if $T$ either has eni-DOP or is eni-deep, then its class
of countable models is Borel complete.  We introduce the
notion of $\lambda$-Borel completeness and prove that
such theories are $\lambda$-Borel complete.  Using this,
we conclude that an $\aleph_0$-stable theory satisfies
$I_{\infty,\aleph_0}(T,\lambda)=2^\lambda$ for all cardinals $\lambda$
if and only if $T$ either has eni-DOP or is eni-deep.
\end{abstract}

\section{Introduction and Preliminaries}

The main theme of the paper will be to produce many disparate models
of an $\aleph_0$-stable theory, assuming some type of non-structure hypothesis.
In all cases, to show the complexity of a model $M$, we concentrate
on the regular types $p\in S(M)$ that have {\em finite dimension in $M$}
i.e., for some (equivalently for every) finite $A\subseteq M$ on which
the regular type is based and stationary, we have $\dim(p|A,M)$ finite.
That is, there is no infinite, $A$-independent set of realizations of $p|A$
in $M$.   Clearly, this notion is  isomorphism invariant.  If $f$ is an isomorphism
between $M$ and $N$, then $p\in S(M)$ has finite dimension in $M$
if and only
if $f(p)$ has finite dimension in $N$.
This yields a criterion for two models to be non-isomorphic:
If two models are isomorphic, their regular types of finite dimension must correspond.

In order to get our non-structure results, we need to identify and analyze both those
regular types that are capable of having finite dimension in a model (which we term {\em eni})
as well as those regular types that are capable of `supporting' an eni type.  Lumped together,
these regular types are called 
{\em eni-active} (see Definition~\ref{lump}) and we call a regular type {\em dull} if it
is not eni-active.  With Proposition~\ref{eninotdull}, we see
that this eni-active/dull partition of regular types
has many equivalent descriptions.  It is particularly useful that this dichotomy
is preserved under the equivalence relation of non-orthogonality.

The paper begins by stating some well-known results about models of $\aleph_0$-stable theories
and then  identifying various species of regular types.
We close Section~1 by proving a structure result for dull types, Proposition~\ref{DULLchain},
that indicates their name is apt.
This result lays the foundation for Theorems~\ref{eniactivetheorem} and \ref{maingap}.

In Section~2, we define the notion of having a `DOP witness' and define many different variants
of `eni-DOP'.  Fortunately, with Theorem~\ref{equiv}, we see that an $\aleph_0$-stable theory
admits one of these variants if and only if it admits them all.  Thereafter, we   
choose the term `eni-DOP' for its brevity.
Theorem~\ref{equiv} also asserts that among $\aleph_0$-stable theories,
eni-DOP is equivalent to the Omitting Types Order Property (OTOP), as well as
to the existence of an independent triple of countable, saturated models over which
the prime model is not saturated.

Our first major result, Theorem~\ref{eniDOPthm}, proves that among
$\aleph_0$-stable theories, those possessing eni-DOP are Borel complete.
The existence of a {\em finite approximation\/} to a DOP witness
(see Subsection~\ref{twotypes}) gives a procedure for constructing
a model $M_G$ to code any bipartite graph $G$.  In such a coding,
the edge set of $G$ corresponds to the types of finite dimension in $M_G$.
However, it is far from obvious how to recover the vertex set of $G$ from
$M_G$.  A weak attempt at this is given in Proposition~\ref{biggy},
where given an isomorphism $f$ between two models $M_G$ and $M_{G'}$,
there is a number $\ell$ (depending largely on $\wt(f(a)/a)$)
so that the image of a complete graph of size $m>\ell$ is almost
complete.   As the number $\ell$ depends on the isomorphism and cannot
be predicted in advance, we obtain our Borel completeness result by
first coding an arbitrary tree $\T$ into a graph $G^*_{\T}$ in which
each node $\eta\in\T$ corresponds to a sequence of finite, complete subgraphs
of arbitrarily large size.  Then, by composing this map with the
coding of graphs into models described above, we obtain a $\lambda$-Borel
embedding of subtrees of $\lambda^{<\omega}$ into models of our theory.
It is noteworthy that had we been able to add finitely many constant symbols to the language,
the proof of Borel completeness in the expanded language would have been much easier.

Once Theorem~\ref{eniDOPthm} has been established, for the remainder
of the paper we assume that $T$ is $\aleph_0$-stable with eni-NDOP.
In Section~5 we introduce and relate several notions of decompositions of a given
model $M$.  
In Definition~\ref{decompdef}, 
decompositions are named [regular, eni, eni-active] according to
the the species of $\tp(a_\nu/M_{\nu^-})$.  With Theorems~\ref{regulartheorem},
\ref{eniactivetheorem}, and \ref{enitheorem} we measure the extent to which
one can recover a model $M$ from a decomposition of it.
Some of these results appear or are implicit in \cite{SHM} and \cite{K},
but are included here to contrast the pros and cons of each species
of decomposition.

In Section 6 we define an $\aleph_0$-stable theory $T$
to be {\em eni-deep} if it has eni-NDOP and some model $M$ has
an eni-active decomposition with an infinite branch.  
With Theorem~\ref{enideepthm}, we prove that 
any $\aleph_0$-stable, eni-deep theory is Borel complete.
The proof uses a major result from \cite{ShL} as a black box.

Finally, in Section 7, we collect our results into Theorem~\ref{maingap},
that characterizes those $\aleph_0$-stable theories that have maximally large
families of $L_{\infty,\aleph_0}$-inequivalent models of every cardinality.

We are grateful to the anonymous referee for mentioning that the class
of eni types need not be closed  under non-orthogonality and for insisting that
the relationship between eni-active types and chains be described more precisely.

\medskip
\centerline{{\bf For the whole of this paper, all theories are $\aleph_0$-stable.}}

\subsection{Preliminary facts about $\aleph_0$-stable theories}

We begin by enumerating several well-known facts about models of
$\aleph_0$-stable theories.

\begin{Definition}  {\em  A non-algebraic, stationary type $p$ is {\em regular} if
$p$ is orthogonal to every forking extension of itself.  $p$ is {\em strongly
regular via $\phi$} if $\phi\in p$ and for every strong type $q$ containing $\phi$,
$q$ is either orthogonal or parallel to $p$.
}
\end{Definition}

It is well known that the binary relation of non-orthogonality is an equivalence relation on the class of stationary,
regular types.  

\begin{Fact}  \label{Fact}
\begin{enumerate}
\item  Over any set $A$, prime and atomic models (indeed, constructible)
models exist and are unique up to isomorphisms over $A$;
\item  If $M$ is a model and $p\not\perp M$, then there is a strongly
regular $q\in S(M)$ non-orthogonal to $p$;
\item  Strongly regular types over models are RK-minimal, i.e., if
$M\preceq N$, $q\in S(M)$ is strongly regular, and there is some
$a\in N\setminus M$ such that $\tp(a/M)\not\perp q$, then
$q$ is realized in $N$;
\item  Any pair $M\preceq N$
of models admits a {\em strongly regular resolution\/} i.e.,
a continuous, elementary chain $\<M_i:i\le\alpha\>$
of elementary substructures of $N$  such
that $M_0=M$, $M_\alpha=N$, and $M_{i+1}$ is prime over $M_i\cup\{a_i\}$, where $\tp(a_i/M_i)$
is strongly regular;
\item  For any complete type $p\in S(M)$ over a model, there is a finite
subset $A\subseteq M$ over which $p$ is based and stationary;
\item  A model is $a$-saturated (i.e., ${\bf F}^a_{\kappa(T)}$-saturated
in the notation of \cite{Shc}) if and only if it is $\aleph_0$-saturated.

\end{enumerate}
\end{Fact}

By combining Fact~\ref{Fact}(2) and (3), we obtain the very useful
`3-model Lemma'.

\begin{Lemma} \label{threemodel}
Suppose $N_0\preceq N_1\preceq M$, $p\in S(N_1)$ is realized
in $M$, and is non-orthogonal to $N_0$.  Then there is a strongly regular
$q\in S(N_0)$
non-orthogonal to $p$ that is realized in $M$ by an element $e$ satisfying
$\fg e {N_0} {N_1}$.
\end{Lemma}

\bp  By Fact~\ref{Fact}(2), choose a strongly regular $q\in S(N_0)$
non-orthogonal to $p$.  Let $q'$ be the non-forking extension of $q$ to
$S(N_1)$.  As $p$ is realized in $M$, it follows from Fact~\ref{Fact}(3)
that $q'$ is realized in $M$ as well.  But any $e$ realizing $q'$ satisfies
$\fg e {N_0} {N_1}$.
\endproof

The following notion is implicit in several proofs of atomicity in \cite{SHM}.

\begin{Definition}  \label{essentiallyfinite} {\em A set $A$ is {\em essentially finite
with respect to a strong type $p$\/} if, for all finite sets $D$ on which $p$ is based
and stationary, there is a finite $A_0\subseteq A$ such that $p|DA_0\vdash p|DA$.
}
\end{Definition}

\begin{Lemma}  \label{basicorth}
Fix a strong type $p$.
If either of the following conditions hold
\begin{enumerate}
\item  $p\perp A$ and $B$ is a (possibly empty) $A$-independent set of finite sets; or
\item  if $A$ is essentially finite with respect to $p$, $p\perp B$, and $\fg A {A\cap B} {B}$
\end{enumerate}
then $A\cup B$ is essentially finite with respect to $p$.
\end{Lemma}

\bp   To establish (1), suppose $B=\{b_i:i\in I\}$ is $A$-independent.
Choose any finite $D$ over which $p$ is based and stationary.  
Now, choose a finite $B_0\subseteq B$ such that $\fg D {AB_0} B$ and then choose a finite $A_0\subseteq A$
such that $\fg {DB_0} {A_0} A$.  We claim that $p|DA_0B_0\vdash p|DAB$.

To see this, first note that since $p\perp A$, we have $p\perp A_0$, 
which coupled with $\fg {DA_0B_0} {A_0} A$
implies $p|DA_0B_0\vdash p|DAB_0$.
But then, since $\fg {DB_0} A {(B\setminus B_0)}$ we obtain 
$p|DAB_0\vdash p|DAB$, proving (1).

To prove (2), write $E:=A\cap B$.  Choose a finite $D$
on which $p$ is based and stationary. 
Choose $B_0\subseteq B$ finite such that $\fg D {B_0A} B$.  As $\fg A E B$
we have
$\fg B {EB_0} {DB_0A}$ and $EB_0\subseteq DB_0A\cap B$.  Choose a finite $A_0\subseteq A$ such that
$\fg {DB_0} {A_0} A$.  Finally, as $A$ is essentially finite with respect to $p$,
choose $A_1\subseteq A$ finite such that  $A_0\subseteq A_1$ and 
$p|DB_0A_1\vdash p|DB_0A$.
Put $D^*:=DB_0A_1$.  As $D^*A=DB_0A$, we have
$\fg B {EB_0} {D^*A}$ and $EB_0\subseteq D^*A$.

To see that $p|D^*\vdash p|DAB$, first note that 
from above, $p|D^*\vdash p|D^*A$.  Also, since $p\perp B$,  $p\perp EB_0$ and 
since $\fg {D^*A} {EB_0} B$ we conclude that
$p|D^*A\vdash p|DAB$.
\endproof

Next, we give a criterion for
$\lambda$-saturation of a model of an $\aleph_0$-stable theory.
For the moment, call a non-algebraic type $p\in S(M)$ {\em $\lambda$-full\/}
if $\dim(p|A,M)\ge\lambda$ for some (every) finite set $A\subseteq M$ on
which $p$ is based and stationary.

\begin{Lemma}  \label{charsat}
For $\lambda$ any infinite
cardinal, a model $M\models T$ is $\lambda$-saturated if and only if every
strongly regular $p\in S(M)$ is $\lambda$-full.
\end{Lemma}

\bp  Left to right is clear, so fix an infinite cardinal $\lambda$
and a model $M$ in which every strongly regular type is $\lambda$-full.
If $M$ is not $\lambda$-saturated, then there is a subset $A\subseteq M$,
$|A|<\lambda$, and a type $q\in S(A)$ that is omitted in $M$.
Among all possible choices, choose $q$ of least Morley rank.
Let $q'\in S(M)$ denote the unique non-forking extension of $q$ to $M$,
let $a$ be any realization of $q'$, and let $N=M[a]$ be any prime model
over $M\cup\{a\}$.  By Fact~\ref{Fact}(2)
 there is an element $b\in N\setminus M$
such that $p=\tp(b/M)$ is strongly regular.  Choose $B\subseteq M$, $|B|<\lambda$,
such that $A\subseteq B$, $p$ is based and stationary over $B$,
and $\tp(a/Bb)$ forks over $B$.  Since $p$ is $\lambda$-full,
there is $b^*\in M$ realizing $p|B$.  Choose any $a^*\in\C$ realizing
$q|B$ with $\tp(a^*/Bb^*)$ forking over $B$.  Now $a^*\not\in M$, lest
$q$ be realized in $M$.  Thus, $r=\tp(a^*/M)$ is non-algebraic, yet
$MR(r)<MR(q)$, hence $r|C$ is realized in $M$ for any $C\supseteq Bb^*$ on
which $r$ is based and stationary and $|C|<\lambda$.  However, any realization
of $r|C$ is a realization of $q$, contradicting $q$ being omitted in $M$.
\endproof

The following Corollary is immediate.

\begin{Corollary}  \label{notcountablesat}
A countable model $M$ is saturated if and only if every strongly regular
$q\in S(M)$ has infinite dimension.
\end{Corollary}

Given two sets $A,B$, we say that {\em $A$ has the Tarski-Vaught property
in $B$,} written $A\subseteq_{TV} B$, if $A\subseteq B$ and every
$L(A)$-formula $\phi(x,a)$ that is realized in $B$ is also realized in $A$.

\begin{Lemma} \label{retain}
\begin{enumerate}
\item  If $B\subseteq_{TV} B'$, then for every $a$, if $\tp(a/B)$ is isolated by the
$L(B)$-formula $\phi(x,b)$, then $\tp(a/B')$ is also isolated by
$\phi(x,b)$.  
\item
 Suppose that $B$ and $C$ are sets with $B$ containing a model $M$
and $\fg B M C$.  Then $B\subseteq_{TV} BC$.  Furthermore, if
$A$ is atomic over $B$, then $\fg {AB} M C$ and $A$ is atomic over $BC$ via the same $L(B)$-formulas.

\item  Suppose that $\<A_i:i<\alpha\>$ and $\<B_i:i<\alpha\>$ are both continuous, increasing subsets of
a model such that each $A_i$ contains and is atomic over $B_i$, and $B_i\subseteq_{TV} B_j$ whenever
$i<j<\alpha$.  Then:
\begin{enumerate}
\item  $B_i\subseteq_{TV} \bigcup B_i$; 
\item  $\bigcup A_i$ is atomic over $\bigcup B_i$; and
\item  If, in addition, each $A_i$ was maximal atomic over $B_i$ inside $N$, then
each $A_i\preceq\bigcup A_i\preceq N$ and $\bigcup A_i$ is maximal atomic over $\bigcup B_i$.
\end{enumerate}
\end{enumerate}
\end{Lemma}

\bp
(1)  is Lemma~XII~1.12(3) of \cite{Shc}, but
we prove it here for convenience.
Let  $\psi(x,b_1,b')$ be any formula over $B'$ with $b_1$ from $B$ and
$b'$ from $B'$.
Let $$\theta(y,z,w):=\forall x\forall x' ([\phi(x,y)\wedge\phi(x',y)]\rightarrow
(\psi(x,z,w)\leftrightarrow\psi(x,z,w)))$$
It suffices to show that $\theta(b,b_1,b')$ holds.
However, if it failed, then  since $b,b_1$ are from $B$
and $B\subseteq_{TV} B'$, we would have
$\neg\theta(b,b_1,b_2)$ for some $b_2$ from $B$.
But this contradicts $\phi(x,b)$ isolating $\tp(a/B)$.

(2)  That $B\subseteq_{TV} BC$ follows from the finite satisfiability of non-forking
over models.  That $\fg {AB} M C$ is a restatement of isolated types being
dominated over models, and the atomicity of $A$ over $BC$ follows from (1).

(3) 
Let $A^*:=\bigcup_{i<\alpha} A_i$ and $B^*:=\bigcup_{i<\alpha} B_i$.
The preservation of the TV-property under continuous chains of sets  is
identical to the preservation of elementarity under continuous chains of models, so $B_i\subseteq_{TV} B^*$ for
each $i$.
To see that $A^*$ is atomic over $B^*$, choose a finite subset $D\subseteq A^*$ and
choose $i<\alpha$ such that $D\subseteq A_i$.  If $\phi(\xbar,b_i)$ isolates $\tp(D/B_i)$,
then by iterating (1), the same formula isolates $\tp(D/B_j)$ for every
$i<j<\alpha$, and hence also isolates 
$\tp(D/B^*)$.

To obtain (c), suppose that each $A_i$ is maximal atomic inside $N$ over
$B_i$.  As there is a prime model $N_i\preceq N$ containing each $A_i$,
the maximality of $A_i$ implies that $A_i\preceq N$, so $A^*\preceq N$ as well. 
To demonstrate that $A^*$ is maximal, choose any $c\in N$ 
such that $A^*c$ is atomic over $B^*$.   Choose $i<\alpha$ such that
both $\tp(c/A^*)$ does not fork and is stationary over $A_i$ and $\tp(c/B_i)$ is isolated.
We will show that $cA_i$ is atomic over $B_i$, which implies $c\in A_i$ by the maximality of $A_i$.
To show this atomicity,
first note that since $A_i$ is atomic over $B_i$ and $B_i\subseteq_{TV} B^*$, it follows from (1)
that $\tp(A_i/B_i)\vdash\tp(A_i/B^*)$, hence
$\fg {A_i} {B_i} {B^*}$.  The transitivity of non-forking implies that
$\fg {cA_i} {B_i} {B^*}$.  Since $cA_i$ is atomic over $B^*$, it 
follows from the Open Mapping Theorem that
$cA_i$ is atomic over $B_i$.
\endproof

Here is an example of a pair of sets with the Tarski-Vaught property.  
It is proved in Lemma~XII~2.3(3) of \cite{Shc}.

\begin{Fact}  \label{Vconfiguration}
Suppose that $M_0,M_1,M_2$ are models with
$\fg {M_1} {M_0} {M_2}$, $N_0$ is a-saturated and
independent from $M_1M_2$ over
$M_0$, $N_1$ is a-prime over $N_0M_1$, and $N_2$ is a-prime over
$N_0M_2$.  Then $M_1M_2\subseteq_{TV} N_1N_2$.
\end{Fact}

\subsection{Species of stationary regular types}

We begin this section by recalling a definition from \cite{ShL}.

\begin{Definition}\label{supporting}
{\em  A stationary, regular type $q$ {\em lies directly above $p$\/}
if there is a non-forking extension $p'\in S(N)$ of $p$ with $N$
$\aleph_0$-saturated, a realization $c$ of $p'$, and an
$\aleph_0$-prime model $N[c]$ over $N\cup\{c\}$ such that
$q\not\perp N[c]$, but $q\perp N$.
A regular type $q$ {\em lies above $p$\/} if there is a sequence
$p_0,\dots,p_n$ of types such that $p_0=p$, $p_n=q$, and
$p_{i+1}$ lies directly above $p_i$ for each $i<n$.  
We say that $p$ {\em supports} $q$ if $q$ lies above $p$.
}
\end{Definition}

The following Lemma gives a sufficient condition for supporting that does not mention $\aleph_0$-saturation.

\begin{Lemma}  \label{notsat}  
Suppose $p\in S(M)$ is regular, $a$ is any realization of $p$, and $M(a)$ is prime over $M\cup\{a\}$.
If a stationary, regular $q\not\perp M(a)$, but $q\perp M$, then $q$ lies directly above $p$, hence
$p$ supports $q$.
\end{Lemma}

\bp  Using Fact~\ref{Fact}(2) and the fact that `lying directly above $p$' is closed under non-orthogonality,
we may assume that $q\in S(M(a))$.  Fix a finite $A\subseteq M(a)$ over which $q$ is based and stationary.
As $M(a)$ is atomic over $M\cup\{a\}$, $\tp(A/Ma)$ is isolated.
Choose any $\aleph_0$-saturated model $N\succeq M$ with $\fg N M a$.  
It follows by Finite Satisfiability that $Ma\subseteq_{TV} Na$, so by Lemma~\ref{retain}(1),
$\tp(A/Na)$ is isolated as well (in fact, by the same formula isolating $\tp(A/Ma)$).
As $\tp(A/Na)$ is $\aleph_0$-isolated, we can choose an $\aleph_0$-prime model $N[a]$ over
$N\cup\{a\}$ that contains $A$.
Now, $p=\tp(a/N)$ is a non-forking extension of $p$
and $q\not\perp N[a]$.  As $A$ is dominated by $a$ over $M$,
$\fg a M N$ and $q\perp M$ we conclude that $q\perp N$.
Thus, $q$ lies directly above $p$, so $p$ supports $q$ by definition.
\endproof

\begin{Definition} \label{eni} {\em  
A stationary, regular type $p$ is {\em eni} (eventually non-isolated) if there is a finite
set $A$ on which $p$ is based and stationary, but $p|A$ is non-isolated.
Such a $p$ is ENI if it is both eni and strongly regular.
}
\end{Definition}

\begin{Definition}  \label{lump} 
{\em  The {\em ENI-active}  types
are the smallest class of stationary, regular types that contain the
ENI types and are closed under automorphisms of the monster model, non-orthogonality,
and supporting.  Similarly,  {\em eni-active} types are the smallest class of stationary, regular
types that are closed under automorphisms, non-orthogonality, and supporting.
}
\end{Definition}

With Proposition~\ref{eninotdull} we will see that every eni type is ENI-active, hence the classes of
ENI-active and eni-active types coincide.
One should note that whereas the class of eni types need not be closed under non-orthogonality,
class of eni-active types is.

\begin{Definition}  {\em  A stationary regular type $p$ is {\em dull\/}
if it is not ENI-active.
}
\end{Definition}

Again, it follows from Proposition~\ref{eninotdull} below that a stationary regular type is dull if and only if it is not eni-active.
Thus, in the notation of Definition~3.7 of \cite{ShL}, if we take  {\bf P}
to be {\em either} the class of ENI types {\em or} the closure of the class of eni types under non-orthogonality, 
then ${\bf P^{\rm active}}$ denotes the class of ENI-active types and ${\bf P^{\rm dull}}$ denotes the dull types.
The remainder of this subsection is aimed at proving Proposition~\ref{eninotdull}.

\begin{Lemma}  \label{infinitedimension}
Suppose that a model $M$ is prime over a finite set $A$, $c$ is a realization of a regular type $p\in S(M)$, 
and $M(c)$ is any prime model over $M\cup\{c\}$.
If $p$ has infinite dimension in $M$, then $M(c)$ is also prime over $A$.
In particular, $M$ and $M(c)$ are isomorphic over $A$.
\end{Lemma}

\bp  First, by increasing $A$ as necessary (while still keeping it finite) we may
assume that $p$ is based and stationary on $A$.
To prove the Lemma, first note that it suffices to find a pair of models
$N\preceq N'\preceq M$ such that $A\subseteq N$ and $N'$ isomorphic over $A$ to
any prime model $N(c)$ over $N\cup\{c\}$.  Indeed, once we have such $N$ and
$N'$, then as they are both countable and atomic over $A$, hence both are isomorphic to
$M$ over $A$.  Thus, $N(c)$ is isomorphic to both $M$ and $M(c)$ over $A$ 
and the Lemma follows.

To produce the submodels $N$ and $N'$, first choose an infinite, $A$-independent
set $J\subseteq M$ of realizations of $p|A$.  Choose a partition $J=J_0\cup J_1$
into disjoint, infinite sets.  Next, choose $B\subseteq M$ to
be maximal subject to the conditions that $AJ_0\subseteq B$ and $\fg B A {J_1}$.
Let $N\preceq M$ be prime over $B$.  

\medskip
\noindent{\bf Claim:}  $N=B$, hence $\fg N A {J_1}$.
\medskip

\bp
Choose any $e\in N$.  As $N$ is atomic over $B$, choose a finite set $C$,
$A\subseteq C\subseteq B$ such that $\tp(e/C)\vdash\tp(e/B)$.
As $J_0\subseteq B$, it follows that $\tp(e/C)\vdash\tp(e/BJ_1)$ [Why?
For $\abar_1$ any tuple from $J_1$, a formula $\phi(x,c,b,\abar_1)\in\tp(e/BJ_1)$ iff there is a cofinite $J_0'\subseteq J_0$ such that $\phi(x,c,b,\abar_0)\in\tp(e/B)$
for some $\abar_0$ from $J_0'$.]  In particular, $\fg {e} B {J_1}$,
which by transitivity implies $\fg {Be} {A} {J_1}$.  Thus, the maximality of $B$
implies that $e\in B$, proving the Claim.
\endproof

Now choose $a\in J_1$ arbitrarily and choose any $N'\preceq M$ to be prime over
$N\cup\{a\}$.  As $\tp(a/A)=\tp(c/A)$ is stationary and both $a$ and $c$ are independent from $N$ over $A$, it follows that $\tp(a/N)=\tp(c/N)$.  In particular,
if we choose $N(c)$ to be any prime model over $N\cup\{c\}$, it will be isomorphic to $N'$ over $N$.  Thus, $N$ and $N'$ are as desired, completing the proof of the Lemma.
\endproof

\begin{Definition}  \label{dullpairdef} 
{\em  A pair of models $M\preceq N$ is a
{\em dull pair} if, for every $d\in N\setminus M$, $\tp(d/M)$ is dull whenever it
is regular.
}
\end{Definition}

\begin{Lemma}  \label{dullpairone}
Suppose $M\preceq N$ is a dull pair, $c\in N\setminus M$ has
$\tp(c/M)$ strongly regular, and $M(c)\preceq N$ is prime over $M\cup\{c\}$.
Then $M(c)\preceq N$ is a dull pair.
\end{Lemma}

\bp  Choose any $d\in N\setminus M(c)$ such that $p:=\tp(d/M(c))$ is regular.
There are two cases.  First, if $p\not\perp M$, then by Lemma~\ref{threemodel}
there is $e\in N\setminus M$ such that $q:=\tp(e/M)$ is strongly regular and
non-orthogonal to $p$.  As $M\preceq N$ is a dull pair, $q$ and hence $p$ must be dull.
On the other hand, suppose that $p\perp M$.  If $p$ were not dull, it would be 
ENI-active, which by Lemma~\ref{notsat} would imply that $\tp(c/M)$ is ENI-active as well,
again contradicting $M\preceq N$ being a dull pair.
\endproof

\begin{Lemma}  \label{dullpairlemma}  Suppose that $M\preceq N$ is a dull pair.
Then for any $M'$ satisfying $M\preceq M'\preceq N$, both $M\preceq M'$ and 
$M'\preceq N$ are dull pairs.
\end{Lemma}

\bp  That $M\preceq M'$ is a dull pair is immediate.
For the other pair, we argue by induction on $\alpha$ that 
\begin{quotation}
For any $M'$ satisfying $M\preceq M'\preceq N$, if there is a strongly regular resolution
$M=M_0\preceq M_1\preceq\dots\preceq M_\alpha=M'$
then $M'\preceq N$ is a dull pair.  
\end{quotation} 
This would suffice by Fact~\ref{Fact}(4), which asserts
the existence of a strongly regular resolution of any $M'$.
When $\alpha=0$ there is nothing to prove.  If $\alpha$ is a non-zero limit ordinal,
then for any $d\in N\setminus M'=M_\alpha$ such that $p:=\tp(d/M_\alpha)$ is regular,
choose $\beta<\alpha$ such that $q:=\tp(d/M_\beta)$ is parallel to $p$.
By induction we have that $q$ is dull, hence $p$ is dull as well.

Finally, assume the inductive hypothesis holds for $\beta$ and 
suppose $M'$ has a strongly regular resolution of length $\alpha=\beta+1$.
By the inductive hypothesis, $M_\beta\preceq N$ is a dull pair, $\tp(c_\beta/M_\beta)$
is strongly regular, and $M'=M_\alpha$ is prime over $M_\beta\cup\{c_\beta\}$.
Thus, Lemma~\ref{dullpairone} implies that $M'\preceq N$ is a dull pair, and our induction is complete.
\endproof

\begin{Proposition}  \label{eninotdull}  
\begin{enumerate}
\item Every eni type is ENI-active;
\item  A type is eni-active if and only if it is ENI-active;
\item  A stationary, regular type is dull if and only if it is not eni-active.
\end{enumerate}
\end{Proposition}

\bp  Once we have proved (1), Clauses (2) and (3) follow immediately from the definitions.  Fix an eni type $p$.  
Choose a finite set $A$ on which $p$ is based and stationary, yet $p|A$ is not isolated.  Let $M$ be prime over $A$.
As $M$ is atomic over $A$, it follows that $M$ omits $p|A$.  Let $e$ be any
realization of $p|M$ and let $M(e)$ be prime over $M\cup\{e\}$.  

By way of contradiction, assume that $p$ were not ENI-active, i.e., $p$ is dull.
We will obtain our contradiction by showing that $M(e)$ is also prime over $A$,
which is a contradiction since $M(e)$ visibly realizes $p|A$.
To obtain this result, we begin with the following Claim.

\medskip\noindent{\bf Claim.}
There is a strongly regular resolution $M=M_0\preceq\dots\preceq M_n=M(e)$ of finite length $n$.
\medskip

\bp  First choose any maximal sequence $M=M_0'\preceq M_1'\preceq M_n'\preceq M(e)$ satisfying the conditions:
(1) $M_{i+1}$ is prime over $M_i\cup\{d_i\}$ where $\tp(d_i/M_i')$ is strongly regular and (2) $\tp(e/M_{i+1}')$ forks
over $M_i'$.  As the sequence of ordinals $\<RM(e/M_i'):i\le n\>$ is strictly decreasing, such a sequence can have at most
finite length.  Also, for any such sequence, we must have $e\in M_n'$, because if not, then by 
Facts~\ref{Fact}(2,3), there would be some strongly
regular type $q\in S(M_n')$ realized in $M(e)\setminus M_n'$ with $q\not\perp \tp(e/M_n')$.  
However, if $d_n$ were any realization of $q$ in $M(e)$ and $M_{n+1}'\preceq M(e)$ were any prime model
over $M_n'\cup\{d_n\}$, we would have $e$ forking with $M_{n+1}'$ over $M_n'$, which would contradict the maximality
of the sequence.  

Thus, any maximal sequence has $e\in M_n'$.  It follows that  $M_n'$ is prime over $M\cup\{e\}$, hence there is an
isomorphism $f:M_n'\rightarrow M(e)$ fixing $M\cup\{e\}$ pointwise.  Then the sequence $\<f(M_i'):i\le n\>$
is a strongly resolution of $M(e)$, completing the proof of the Claim.
\endproof

Fix such a strongly regular resolution
$M=M_0\preceq M_1\preceq M_n=M(e)$ where $M_{i+1}$ is prime over $M_i\cup\{c_i\}$ and $\tp(c_i/M_i)$ is strongly
regular.  Next, note that $e$ forks over $M$ with any $d\in M(e)\setminus M$, so if $\tp(d/M)$ is regular it must
be non-orthogonal to $p$ and hence dull.  That is, $M\preceq M(e)$ is a dull pair.  It follows from Lemma~\ref{dullpairlemma}
that $M_i\preceq M(e)$ is also a dull pair whenever $i<n$.  

Using this, we complete the proof by showing, by induction on $i\le n$, that each $M_i$ is prime over $A$.
When $i=0$ this is immediate by hypothesis.  So fix $i<n$ and assume that $M_i$ is prime over $A$.
Let $q_i:=\tp(c_i/M_i)$.  Choose a finite set $B$, $A\subseteq B\subseteq M_i$ on which $q_i$ is based and
stationary.  Note that $M_i$ is prime over $B$ as well.  
As $M_i\preceq M(e)$ is a dull pair, $q_i$ is strongly regular and dull.  In particular, $q_i$
is not ENI, hence $q_i$ has infinite dimension in $M_i$.  
Thus, $M_{i+1}$ is prime over $B$ by Lemma~\ref{infinitedimension}.
However, as $\tp(B/A)$ is isolated, it follows that $M_{i+1}$ is prime over $A$ as well.
\endproof

\subsection{On dull types}

We begin by defining a strong notion of substructure.

\begin{Definition} {\em  $N$ is an {\em $L_{\infty,\aleph_0}$-substructure of $M$\/} 
if $N\preceq M$ and for all finite $A\subseteq N$,
$$(N,a)_{a\in A}\bfequiv (M,a)_{a\in A}$$
}
\end{Definition}

The paradigm of this notion is when $M$ is atomic and $N\preceq M$.  In this case,
both $N$ and $M$ are atomic, hence back-and-forth equivalent, over every finite $A\subseteq N$.
With Proposition~\ref{DULLchain}, we prove that this stronger notion of substructure holds
for every dull pair $N\preceq M$.
We begin with a Lemma which gets its strength when coupled with
Lemma~\ref{basicorth}.

\begin{Lemma} \label{atomicextension}
Suppose $A\subseteq C$ is essentially finite with respect to a regular, stationary,
but not eni type $p\in S(C)$.  If $C$ is atomic over $A$, then so is $C\cup\{e\}$
for any realization $e$ of $p$.  
\end{Lemma}

\bp  It suffices to show that $De$ is atomic over $A$ for every finite $D\subseteq C$.  So fix any
finite $D\subseteq C$.  Choose a finite $D^*$, $D\subseteq D^*\subseteq C$
with  $p$
 based and stationary on $D^*$.
As $A$ is essentially finite with respect to $p$, choose a finite $A_0\subseteq A$
such that $p|D^*A_0\vdash p|D^*A$.
Since $p$ is not eni and $D^*A_0$ is finite, $\tp(e/D^*A_0)$ is isolated.  
Coupled with the fact that $D^*$ is atomic over $A$, this implies
that $D^*e$ and hence $De$ is atomic over $A$, as required.
\endproof

\begin{Lemma}  \label{dulltechnical}  
Suppose that $N\preceq N(c)$, where $N(c)$ is prime over $N\cup\{c\}$ and $c$ realizes a dull type $p\in S(N)$.
Then for every finite set $A$, there is $M\preceq N$ over which $p$ is based and an infinite Morley sequence $J\subseteq N$
in $p|M$ such that
\begin{itemize}
\item  $A\subseteq M$;
\item  $N$ is atomic over $M\cup J$; and
\item  $N(c)$ is atomic over $M\cup J\cup\{c\}$.
\end{itemize}
\end{Lemma}  

\bp
Without loss, we may assume that $p$ is based and stationary on $A$.
As $p$ is not eni, $p$ has infinite dimension in $N$, so we can find 
an infinite Morley sequence $J^*\subseteq N$ in $p|A$.  Partition $J^*$ into two 
disjoint, infinite
pieces $J^*=J_0\cup J$.
Arguing as in the proof of Lemma~\ref{infinitedimension}, 
choose $B\subseteq N$ maximal subject to the conditions
(1) $AJ_0\subseteq B$ and (2) $\fg B {A} J$.
Just as in \ref{infinitedimension}, 
$B$ is the universe of an elementary substructure, which we denote as $M$.
Clearly, $A\subseteq M$.

\medskip\noindent{\bf Claim 1.}  $M\preceq N$ is a dull pair.
\medskip

\bp  Choose any $e\in N$ such that $q=\tp(e/M)$ is regular.  We show that any such $q$ must be non-orthogonal to $p$
and hence be dull.  If this were not the case, then we would have $\fg e M J$, which would contradict the maximality of
$M$.
\endproof

\medskip\noindent{\bf Claim 2.}  $N$ is atomic over $M\cup J$.
\medskip

\bp  Choose $N_0\preceq N$ to be maximal atomic over $M\cup J$.  We argue that $N_0=N$.
If this were not the case, choose $e\in N$ such that $q:=\tp(e/N_0)$ were regular.  As $M\preceq N$ is a dull
pair, it follows from Lemma~\ref{dullpairlemma} that $q$ is dull and hence not eni by 
Proposition~\ref{eninotdull}.  We argue by cases.  First, if $q$ were non-orthogonal to $M$,
then by Lemma~\ref{threemodel} there would be $d\in N\setminus M$ such that $\fg d M {N_0}$
which, since $J\subseteq N_0$, would contradict the maximality of $M$.  
On the other hand, if $q\perp M$, then by Lemmas~\ref{basicorth}(1) and 
\ref{atomicextension}
we would have $N_0\cup\{e\}$ atomic over $M\cup J$, which  contradicts the maximality of $N_0$.

\medskip\noindent{\bf Claim 3.}  $M$ is maximal in $N(c)$ such that $\fg M A {Jc}$.
\medskip

\bp  First, it is clear that $\fg M A {Jc}$ by the defining property of $M$ and because $\fg c A N$.
The verification of the maximality of $M$ inside $N(c)$ is an exercise in non-forking.
Namely, choose any $e\in N(c)$ such that $\fg {eM} A {Jc}$.  
As $J\cup\{c\}$ is independent over $A$, we have 
$\fg {eMc} A {J}$, hence $\fg {ec} M {J}$.  As $N$ is atomic over
$M\cup J$ by Claim~2, we obtain $\fg {ec} M N$.
Combining this with the fact that $\fg e M c$ yields $\fg e M {Nc}$,
hence $\fg e M {N(c)}$.  Thus, $e\in M$ as required.
\endproof

We finish by using analogues of the proofs of Claims 1 and 2 (using $Jc$ in place of $J$)
to prove that
$M\preceq N(c)$ is a dull pair and that $N(c)$ is atomic over $MJc$.
\endproof

\begin{Lemma} \label{iso}
Suppose that
$N\models T$,  $\tp(c/N)$ is dull, and $N(c)$ is any prime model over $N\cup\{c\}$.
Then $N$ is an $L_{\infty,\aleph_0}$-elementary substructure of $N(c)$.
\end{Lemma}

\bp  Given $N\preceq N(c)$ and a finite $A$, by enlarging $A$ slightly we may
assume that $p=\tp(c/N)$ is based and stationary on $A$.
Apply Lemma~\ref{dulltechnical} to obtain $M\preceq N$ and $J$ such that
$A\subseteq M$, $N$ is atomic over $MJ$, and $N(c)$ is atomic over $MJc$.

Let $g:MJ\rightarrow MJc$ be any elementary bijection that is the identity on $M$.

But now, we show that $(N,a)_{a\in M}\bfequiv (N(c),a)_{a\in M}$
by exhibiting the back-and-forth system 
\begin{quotation}
${\cal F}=\{$all
finite partial functions $f:N\rightarrow N(c)$ such that
$f\cup g$ is elementary$\}$
\end{quotation}
The verification that $\F$ is a back-and-forth system is akin to the verification
that any two atomic models of a complete theory are back-and-forth equivalent.
\endproof

The following Proposition follows by iterating Lemma~\ref{iso}:

\begin{Proposition} \label{DULLchain}
Suppose that $N\preceq M$ is a dull pair.
Then $N$
is an
 $L_{\infty,\aleph_0}$-elementary
substructure of $M$.
\end{Proposition}

\bp  Choose a strongly regular resolution 
$N=N_0\preceq N_1\preceq\dots\preceq N_\alpha=M$
As $\tp(c_{i+1}/N_i)$ is dull for each $i<\alpha$, it follows from Lemma~\ref{iso}
that $N_i$
is an
 $L_{\infty,\aleph_0}$-elementary
substructure of $N_{i+1}$.
\endproof

\section{eni-DOP and equivalent notions} \label{eniDOPsection}

We begin with a central notion of \cite{ShL} and contrast it with a slight strengthening.

\begin{Definition} {\em  A stationary, regular type $p$ has a {\em DOP witness\/}
if there is a quadruple $(M_0,M_1,M_2,M_3)$ of models,
where $(M_0,M_1,M_2)$ form an independent triple of $a$-models,
$M_3$ is $a$-prime over $M_1\cup M_2$, $p$ is based on $M_3$, but
$p\perp M_1$ and $p\perp M_2$.
A {\em prime DOP witness\/} for $p$ is the same, except that we require that
$M_3$ be prime over $M_1\cup M_2$ (as opposed to a-prime).
}
\end{Definition}

Visibly, among stationary, regular types,
having either a DOP witness or a prime DOP witness
is invariant under parallelism and automorphisms of the monster model $\C$.

Recall that by Fact~\ref{Fact}(6), an a-model is simply an
$\aleph_0$-saturated model.  As in \cite{ShL}, we are free to vary the
amount of saturation of the models $(M_0,M_1,M_2)$.  

\begin{Lemma}  \label{DOPwitness}
The following are equivalent for a stationary regular type $p$.
\begin{enumerate}
\item  $p$ has a prime DOP witness;
\item  There is a quadruple $(M_0,M_1,M_2,M_3)$ of models such
that $(M_0,M_1,M_2)$ form an independent triple, $M_3$ is prime over
$M_1\cup M_2$, $p$ based on $M_3$, but $p\perp M_1$ and $p\perp M_2$;
\item  Same as (2), but with
$\dim(M_1/M_0)$ and $\dim(M_2/M_0)$  both finite;
\item Same as (2), but with $\dim(M_1/M_0)=\dim(M_2/M_0)=1$;
\item $p$ has a prime DOP witness $(M_0,M_1,M_2,M_3)$  with
$\dim(M_1/M_0)=\dim(M_2/M_0)=1$.
\end{enumerate}
\end{Lemma}

\bp  $(1)\Rightarrow(2)$ is immediate.

$(2)\Rightarrow(3)$:  Let $(M_0,M_1,M_2,M_3)$ be any witness to (2).
Choose a finite $d\subseteq M_3$ over which $p$ is based and stationary.
As $M_3$ is prime over $M_1\cup M_2$, choose finite $b\subseteq M_1$ and
$c\subseteq M_2$ such that $\tp(d/M_1M_2)$ is isolated by a formula $\phi(x,b,c)$.
Let  $M_1'\preceq M_1$ be prime over $M_0b$, let
$M_2'\preceq M_2$ be prime over $M_0c$, and let $M_3'\preceq M_3$
be prime over $M_1'\cup M_2'$ with $d\subseteq M_3'$.  
Then $(M_0,M_1',M_2',M_3')$
are as required in (3).

$(3)\Rightarrow(4)$:  Assume that (3) holds.  Among all possible
quadruples of models witnessing (3), choose a triple $(M_0,M_1,M_2,M_3)$
with $\dim(M_1/M_0)+\dim(M_2/M_0)$ as small as possible.
Clearly, we cannot have either $M_0=M_1$ or $M_0=M_2$, so
$\dim(M_1/M_0)$ and $\dim(M_2/M_0)$ are each at least one.
We argue that the minimum sum occurs when
$\dim(M_1/M_0)=\dim(M_2/M_0)=1$.  Assume this were not the case.
Without loss, assume that $\dim(M_1/M_0)\ge 2$.  Choose an element
$e\in M_1\setminus M_0$ such that $\tp(e/M_0)$ is strongly regular
and let $M_1'\preceq M_1$ be prime over $M_0\cup\{e\}$.
Let $M_3'\preceq M_3$ be prime over $M_1'\cup M_2$.
There are two cases.  On one hand, if $p\not\perp M_3'$,
then by e.g., Claim~X.1.4 of \cite{Shc},
choose an automorphic copy $p'$ of $p$ that is based on $M_3'$
with $p\not\perp p'$.
Then an automorphic copy of the quadruple $(M_0,M_1',M_2,M_3')$ 
contradicts the minimality of
our choice.  On the other hand, if $p\perp M_3'$, then the quadruple
$(M_1',M_1,M_3',M_3)$ directly contradicts the minimality of our choice.

$(4)\Rightarrow(5)$:  Let $(M_0,M_1,M_2,M_3)$ be any witness to (4).
Choose a finite $d\subseteq M_3$ over which $p$ is based and stationary.
Now, choose an a-model $N_0$ satisfying $\fg {N_0} {M_0} {M_3}$,
and choose a-prime models $N_1$ and $N_2$ over $N_0\cup M_1$ and
$N_0\cup M_2$, respectively.  As $\fg {M_3} {M_\ell} {N_\ell}$ for both
$\ell=1,2$, it follows that $p\perp N_1$ and $p\perp N_2$.
Also, as $M_1M_2\subseteq_{TV} N_1N_2$ by Fact~\ref{Vconfiguration},
it follows 
Lemma~\ref{retain}(1)  
that $\tp(d/N_1N_2)$ is isolated.
Choose a prime model $N_3$ over
$N_1\cup N_2$
that contains $d$.  Then $(N_0,N_1,N_2,N_3)$  is a prime DOP witness
for $p$ with $\dim(N_1/N_0)=\dim(N_2/N_0)=1$.

$(5)\Rightarrow(1)$ is immediate.
\endproof

\begin{Definition} \label{eni-DOP}
{\em  $T$ has {\em eni-DOP\/}
if some eni type
$p$ has a DOP witness.  Similarly,
$T$ has {\em ENI-DOP} (respectively, {\em eni-active DOP\/}) if some
ENI-type (respectively, eni-active type) has a DOP witness.
}
\end{Definition}

It is fortunate, at least for
the exposition, that $T$ having any of the three preceding notions
are equivalent.  
In fact, this equivalence extends much further.
Recall that a stable theory has the {\em Omitting Types Order Property (OTOP)}
if there is a type $p(x,y,z)$ (where $x,y,z$ denote finite tuples of variables)
such that for any cardinal $\kappa$ there is a model $M^*$ and
a sequence $\<(b_\alpha,c_\alpha):\alpha<\kappa\>$ such that for all
$\alpha,\beta<\kappa$,
$$ M^*\ \hbox{realizes}\ p(x,b_\alpha,c_\beta) \quad\hbox{if and only if}
\quad \alpha<\beta$$

\begin{Theorem}  \label{equiv}
The following are equivalent for an $\aleph_0$-stable
theory $T$:
\begin{enumerate}
\item $T$ has eni-DOP;
\item $T$ has ENI-DOP;
\item $T$ has eni-active DOP;
\item There is an independent triple $(M_0,M_1,M_2)$ of countable,
saturated models such that some (equivalently every) prime model
over $M_1\cup M_2$ is not saturated;
\item There is an independent triple $(N_0,N_1,N_2)$ of countable saturated models
and strongly regular types $p,q\in S(N_0)$
such that
$N_1$ is $\aleph_0$-prime over $N_0$ and a realization $b$ of $p$,
$N_2$ is $\aleph_0$-prime over $N_0$ and a realization $c$ of $q$,
and if $N_3$ is prime over $N_1N_2$, then there is a finite $d$ satisfying $\{b,c\}\subseteq d\subseteq N_3$
and  an ENI type $r(x,d)$ that is omitted in $N_3$
and orthogonal to both
$N_1$ and $N_2$;
\item $T$ has OTOP.
\end{enumerate}
\end{Theorem}

\bp  If we let {\bf P} denote any of eni, ENI, or eni-active, then
it follows from Proposition~\ref{eninotdull} that ${\bf P^{\rm active}}$
(which is the closure of {\bf P} under automorphisms, non-orthogonality and `supporting'
within the class of stationary, regular types)
would be the set of eni-active types.
Thus, Clauses (1), (2) and (3) are equivalent by way of Corollary~3.9 of \cite{ShL}.

$(1)\Rightarrow (4)$:
Suppose that $(M_0,M_1,M_2,M_3)$ is a DOP witness for an eni type $p$
with each model countable and saturated.
Let $N\preceq M_3$ be prime over $M_1\cup M_2$, and by way of contradiction,
assume that $N$ is saturated.  Then as $N$ and $M_3$ are isomorphic
over $M_1\cup M_2$,
by replacing  $p$ by a conjugate type,
we may  assume that $p\in S(N)$.
We will contradict the saturation of $N$ by
finding a finite subset $D^*\subseteq N$ on which $p$ is based and stationary,
but $p|D^*$ is omitted in $N$.

First, since $p$ is eni and $N$ is saturated, choose a finite $D\subseteq N$
on which $p$ is based and stationary, but $p|D$ is not isolated.
As $p\perp M_1$ and $p\perp M_2$, it follows from Lemma~\ref{basicorth} that $M_1M_2$
is essentially finite with respect to $p$.  Thus, there is a finite $D^*\subseteq DM_1M_2$
containing $D$
such that $p|D^*\vdash p|DM_1M_2$.  As $p|D^*$ is a non-forking extension of $p|D$, it cannot
be isolated.
We argue that $p|D^*$ cannot be realized in $N$.  Suppose $c\in N$ realized $p|D^*$.  
Then, as $cD^*$ is atomic over $M_1M_2$, we would have $\tp(c/D^*M_1M_2)$ isolated.
However, since $\tp(c/D^*)\vdash\tp(c/D^*M_1M_2)$, we have $\fg c {D^*} {M_1M_2.}$
Thus, the Open Mapping Theorem would imply that $\tp(c/D^*)$ is isolated, which is a contradiction.

$(4)\Rightarrow (5)$:
Let $(M_0,M_1,M_2)$ exemplify (4), and fix a prime model $M_3$ over
$M_1\cup M_2$.  As $M_3$ is not saturated, by Lemma~\ref{charsat}
there is an ENI $r\in S(M_3)$
of finite dimension in $M_3$.

\medskip
\noindent{\bf Claim.}  $r$ is orthogonal to both $M_1$ and $M_2$.

\bp  As the cases are symmetric, assume by way of contradiction that
$r\not\perp M_1$.
 By Fact~\ref{Fact}(2) there is a strongly regular
$p\in S(M_1)$ nonorthogonal to $r$.  Choose a finite $A\subseteq M_3$ such that
$r$ is based, stationary and strongly regular over $A$, and $r|A$ is omitted in $M_3$.
Choose a finite $B\subseteq M_1$ over which $p$ is based, stationary
and strongly regular, and let
$r'$ and $p'$ be the unique nonforking extensions of $r|A$ and $p|B$ to $AB$.
Since $M_1$ is saturated, $\dim(p|B,M_1)$ is infinite, hence $\dim(p',M_3)$ is infinite
as well.  Thus, $\dim(r',M_3)$ is also infinite, contradicting the fact that $r|A$ is omitted in
$M_3$.  \endproof

Thus, $r$ has a prime DOP witness by Lemma~\ref{DOPwitness}(2).
But now, Lemma~\ref{DOPwitness}(5) gives us the configuration we need.

$(5)\Rightarrow(2)$:  Given the triple $(M_0,M_1,M_2)$ and the type $r$ in (5),
choose an a-prime model $M_3$ over $M_1\cup M_2$.
Then $(M_0,M_1,M_2,M_3)$
is a DOP witness for the ENI type $r$.

$(5)\Rightarrow(6)$:  Given the data from (5),
let $w(x,u,y,z)$ be the type asserting that $y$ and $z$ are $M_0$-independent
solutions of $p$ and $q$, respectively, $\phi(u,y,z)$ isolates $\tp(d/M_1M_2)$
and $r(x,u)$.
We argue that the type $\exists u w(x,u,y,z)$ witnesses OTOP.
To see this, fix any cardinal $\kappa$.  Choose $\{b_i:i<\kappa\}\cup\{c_j:j<\kappa\}$ to be $M_0$-independent, where
$\tp(b_i/M_0)=p$ and $\tp(c_j/M_0)=q$ for all $i,j\in\kappa$.
For each $i,j$, let $M_1(b_i)$ be prime over $M_0\cup\{b_i\}$ and
$M_2(c_j)$ be prime over $M_0\cup\{c_j\}$, and let $\Mbar$
be prime over the union of these models.
Now, for each pair $(i,j)$, choose a witness $d_{i,j}$ to $\phi(u,b_i,c_j)$
from $\Mbar$ and let $r_{i,j}$ be shorthand for $r(x,d_{i,j})$.
It is easily checked that all of the types $r_{i,j}$ are orthogonal.

For each pair $(i,j)$ with $i\le j$, choose a realization $e_{i,j}$ of
$r_{i,j}$, and let $M^*$ be prime over $\Mbar\cup\{e_{i,j}:i\le j<\kappa\}$.
Then, because of the orthogonality of the $r_{i,j}$, $M^*$ realizes
$\exists u w(x,u,b_i,c_j)$ if and only if $i\le j$.

$(6)\Rightarrow(1)$:  This is Corollary~\ref{notop}.  (There is no circularity.)
\endproof

\section{$\lambda$-Borel completeness}
Throughout this section, we {\bf fix a cardinal $\lambda\ge\aleph_0$}.
We consider only models of size $\lambda$, typically those whose universe
is the ordinal $\lambda$, in a language of size $\kappa\le\lambda$.
{\em For notational simplicity, we only consider relational languages.}
Although it would be of interest to explore this notion in more generality,
here we only study classes ${\bf K}$ of $L$-structures that are closed under
$\bfequiv$ and study the complexity of ${\bf K}/\bfequiv$.

\begin{Definition} {\em   For any (relational) language $L$ with at most $\lambda$
symbols, let $\Lplusminus:=L\cup\{\neg R:R\in L\}$,
and let $S^\lambda_L$ denote the set of $L$-structures $M$ with universe
$\lambda$.  Let
$$L(\lambda):=\{R(\alphabar):R\in\Lplusminus,\ \alphabar\in {^{{\rm arity(R)}}\lambda}\}$$
and endow $S^\lambda_L$ with the topology formed by
letting
$$\B:=\{U_{R(\alphabar)}:R(\alphabar)\in L(\lambda)\}$$ be a subbasis,
where $U_{R(\alphabar)}=\{M\in S^\lambda_L:M\models R(\alphabar)\}$.
}
\end{Definition}

\begin{Definition} {\em  Given a language $L$ of size at most
$\lambda$, a set $K\subseteq S^\lambda_L$ is {\em $\lambda$-Borel}
if, there is a $\lambda$-Boolean combination $\Psi$ of
$L(\lambda)$-sentences (i.e., a propositional
$L_{\lambda^+,\aleph_0}$-sentence of $L(\lambda)$) such that
$$K=\{M\in S^\lambda_L:M\models\Psi\}$$
Given two languages $L_1$ and $L_2$, a function
$f:S^\lambda_{L_1}\rightarrow S^\lambda_{L_2}$ is
{\em $\lambda$-Borel} if the inverse image of every (basic) open set
is $\lambda$-Borel.
}
\end{Definition}

That is,  $f:S^\lambda_{L_1}\rightarrow S^\lambda_{L_2}$ is
$\lambda$-Borel if and only if for every $R\in L_2$ and
$\betabar\in {^{\rm arity(R)}\lambda}$, there is a $\lambda$-Boolean
combination $\Psi_{R(\betabar)}$ of $L_1(\lambda)$-sentences such that
for every $M\in S^\lambda_{L_1}$,
$f(M)\models R(\betabar)$ if and only if $M\models \Psi_{R(\betabar)}$.

As two countable structures are isomorphic if and only if
they are $\bfequiv$, a moment's thought tells us that when 
$\lambda=\aleph_0$, the notions
of  $\aleph_0$-Borel sets and functions defined above are equivalent to the
usual notion of  Borel sets and functions.

\begin{Definition}  {\em Suppose that $L_1,L_2$ are relational languages 
with at most
$\lambda$ symbols, and for $\ell=1,2$, $K_\ell$ is a $\lambda$-Borel
subset of $S^\lambda_{L_\ell}$ that is invariant under $\bfequiv$.
 We say that {\em $(K_1,\bfequiv)$ is
$\lambda$-Borel reducible to $(K_2,\bfequiv)$}, written
$$(K_1,\bfequiv)\le_\lambda^B (K_2,\bfequiv)$$
if there is a $\lambda$-Borel function $f:S^\lambda_{L_1}\rightarrow
S^\lambda_{L_2}$ such that $f(K_1)\subseteq K_2$ and, for all
$M,N\in K_1$,
$$M\bfequiv N\qquad\hbox{if and only if}\qquad f(M)\bfequiv f(N)$$
}
\end{Definition}

\begin{Definition}  {\em  A class $K$ is {\em $\lambda$-Borel complete 
for $\bfequiv$\/} if $(K,\bfequiv)$ is a maximum with respect to
$\le^B_\lambda$.  We call a theory $T$ $\lambda$-Borel complete for $\bfequiv$
if $Mod_\lambda(T)$, the class of models of $T$ with universe $\lambda$, is
$\lambda$-Borel complete for $\bfequiv$.
}
\end{Definition}

To illustrate this notion, we prove a series of Lemmas, culminating in a generalization
of Friedman and Stanley's \cite{FS} result that subtrees of $\omega^{<\omega}$ are Borel complete.
We make heavy use of the following characterizations of $\bfequiv$-equivalence of
structures of size $\lambda$.

\begin{Fact} \label{Levy}
If $|L|\le\lambda$, the following conditions are equivalent for
 $L$-structures $M$ and $N$ that are both of size $\lambda$.
\begin{enumerate}
\item  $M\bfequiv N$;
\item  $M$ and $N$ satisfy the same $L_{\lambda^+,\aleph_0}$-sentences;
\item  If $G$ is a generic filter of the Levy collapsing poset $Lev(\aleph_0,\lambda)$, then in $V[G]$ there is an isomorphism $h:M\rightarrow N$ of
countable structures.
\end{enumerate}
\end{Fact}

For all $\aleph_0\le\kappa\le\lambda$, let $L_\kappa$ be the language consisting
of the binary relation $\trianglelefteq$ and $\kappa$ unary predicate symbols
$P_i(x)$.  Let
 $\kappa\,CT_\lambda$  denote the class of all $L_\kappa$-trees with universe
$\lambda^{<\omega}$, colored by the predicates $P_i$.

\begin{Lemma} \label{leftpart} For any (relational)
language $L$ satisfying $|L|\le\kappa\le\lambda$,
$$(S^\lambda_{L},\bfequiv)\le^B_\lambda (\kappa\,CT_\lambda,
\bfequiv)$$
\end{Lemma}

\bp  For each $n\in\omega$, let $\<\phi_{n,i}(\xbar):i<\gamma(n)\le\kappa\>$
be a maximal set of pairwise non-equivalent quantifier-free $L$-formulas
with $\lg(\xbar)=n$.  As well, fix a bijection
$\Phi:\omega\times\kappa\rightarrow\kappa$.

Now, given any $L$-structure $M\in S^\lambda_L$, first note that since the universe of $M$ is $\lambda$, the finite sequences from $M$ naturally form
a tree isomorphic to $\lambda^{<\omega}$ under initial segment.

So $f(M)$ will consist of this tree, with $\trianglelefteq$ interpreted as the
initial segment relation.  Furthermore, for each $j\in\kappa$, choose
$(n,i)\in\omega\times\kappa$ such that $\Phi(n,i)=j$.  If $i<\gamma(n)$,
then put
$$P_j^{f(M)}:=\{\alphabar\in\lambda^n:M\models\phi_{n,i}(\alphabar)\}$$
(if $i\ge\gamma(n)$, then for definiteness, say that $P_j$ always fails on $f(M)$).

Choose any $M,N\in S^\lambda_L$.
It is apparent from the construction that if $M\bfequiv N$, then $f(M)\bfequiv f(N)$.
The other direction is more interesting.  Suppose that $f(M)\bfequiv f(N)$.
Consider the Levy collapsing forcing, $Lev(\aleph_0,\lambda)$, that, for any
generic filter $G$, $V[G]$ includes a bijection $g:\omega\rightarrow\lambda$.
We work in $V[G]$.  Note that both $f(M)$ and $f(N)$ are
$\bfequiv$-equivalent, countable structures.  Thus, in $V[G]$,
fix an $L_\kappa$-isomorphism $h:f(M)\rightarrow f(N)$.
Using $h$, in $\omega$ steps we construct two branches
$\eta,\nu\in \lambda^\omega$, where we think of $\eta$ as a branch through
$f(M)$, while $\nu$ is a branch through $f(N)$, satisfying the following
three conditions:
\begin{itemize}
\item  For each $n\in\omega$, $h(\eta(n))=\nu(n)$;
\item $\{g(n):n\in\omega\}\subseteq\dom(\eta)$; and
\item $\{g(n):n\in\omega\}\subseteq\dom(\nu)$.
\end{itemize}
Let $F=\{(\eta(n),\nu(n)):n\in\omega\}$.
As $\{g(n):n\in\omega\}$ is all of $\lambda$, it follows
that $\dom(F)=\lambda$ and $\range(F)=\lambda$.
Furthermore, since $h(\eta(n))=\nu(n)$, it follows that
$P_j(\eta(n))\leftrightarrow P_j(\nu(n))$ for each $j$.
Thus, for each $n$, the $L$-quantifier free types of
$\<\eta(i):i<n\>$ and $\<\nu(i):i<n\>$ are the same.
In particular, it follows that $F$ is a bijection from $\lambda$ to $\lambda$
that preserves $L$-quantifier-free types.  Thus, $F:M\rightarrow N$ is
an isomorphism.

Of course, the isomorphism $F\in V[G]$, but it follows easily by absoluteness
that $M\bfequiv N$ in $V$.
\endproof

\begin{Definition} {\em
 Given any trees $T$ and $\{S_\eta:\eta\in T\}$,
we form the tree $T^*(S_\eta:\eta\in T)$ that `attaches $S_\eta$ to $T$ at $\eta$'
as follows:

The universe of $T^*(S_\eta:\eta\in T)$ (which, for simplicity,
 we write as $T^*$ below) is the disjoint union
of
$$T\sqcup\bigsqcup_{\eta\in T} S_\eta\setminus\{\<\>\}$$
and, for $u,v\in T^*$, we say $u\le_{T^*} v$ if and only if one
of the following clauses hold:
\begin{itemize}
\item  $u,v\in T$ and $u\trianglelefteq_T v$; or
\item  for some $\eta\in T$, $u,v\in\S_\eta\setminus\{\<\>\}$ and
$u\trianglelefteq_{S_\eta} v$; or
\item  $u\in T$, $v\in S_\eta\setminus\{\<\>\}$ and $u\trianglelefteq_T\eta$.
\end{itemize}
}
\end{Definition}

Note that in particular, elements from distinct $S_\eta$'s are incomparable, and
that no element of any $S_\eta$ is `below" any element of $T$.
It is easily checked that if $T$ and each of the $S_\eta$'s are subtrees of
$\lambda^{<\omega}$, then the attaching tree $T^*(S_\eta:\eta\in T)$
can also be construed as being a subtree of $\lambda^{<\omega}$.

\begin{Definition} {\em
 A {\em subtree of $\lambda^{<\omega}$} is simply
a non-empty subset of $\lambda^{<\omega}$ that is closed under initial
segments.  Given a subtree $T$ of $\lambda^{<\omega}$, an element
$\eta\in T$  is {\em contained in a branch} if there is some $\nu\in\lambda^\omega$ extending $\eta$ such that $\nu(n)\in T$ for
every $n\in\omega$.  A subtree $T$ of $\lambda^{<\omega}$ is {\em special}
if, for every $\eta\in T$ that is contained in a branch, $\eta$ has no
immediate successors that are leaves (i.e., every immediate successor
of $\eta$ has a successor in $T$).
}
\end{Definition}

\begin{Lemma} \label{rightpart}
$(\aleph_0\, CT_\lambda,\bfequiv)\le^B_\lambda\,
(\hbox{Special subtrees of $\lambda^{<\omega}$},
\bfequiv)$

\end{Lemma}

\bp  Fix a bijection $\Phi:\omega\times\omega\rightarrow\omega\setminus\{0,1\}$.
Let $T_0$ be the tree $\lambda^{<\omega}$.

Also, given any subset $V\subseteq\omega$, let $S_V$ be the rooted
tree consisting of one copy of the tree $\omega^{\le m}$ for
each $m\in V$.  Other than being joined at the root, the copies of
$\omega^{\le m}$ are disjoint.

Now, suppose we are given $M\in\aleph_0\, CT_\lambda$, i.e.,
the tree $(\lambda^{<\omega},\trianglelefteq)$, adjoined by
countably many unary predicates $P_j(x)$.
We construct a special tree $f(M)$ as follows:

First, form the tree $T_0=\lambda^{<\omega}$.
For each $\eta\in T_0$, let
$$V(\eta):=\{\Phi(n,j):M\models P_j(\eta)\}$$
where $n=\lg(\eta)$.
Note that each $V(\eta)\subseteq\omega\setminus\{0,1\}$.
Let $f(M)$ be the tree $T_0(S_{V(\eta)}:\eta\in T_0)$.
By the remark above, as each of $T$, $T_0$ and each $S_V$ is a subtree
of $\lambda^{<\omega}$,  $f(M)$ is also a subtree of $\lambda^{<\omega}$.
Furthermore, note that $T_0$ is recognizable in $f(M)$
as being precisely those elements of $f(M)$ that are contained in an infinite
branch.
Moreover, for every element $\eta\in f(M)$ that is not contained in an infinite
branch,
there is a uniform bound on the lengths of $\nu\in f(M)$ extending $\eta$.
Combining this with the fact that $1\not\in V(\eta)$ for any
$(\eta)\in T_0$, we conclude that $f(M)$ is special.

It is easily verified by the construction that if $M\bfequiv N$, then
$f(M)\bfequiv f(N)$.
Conversely, suppose that $M,N\in\aleph_0\, CT_\lambda$ and that
$f(M)\bfequiv f(N)$.
Choose any generic filter $G$ for the Levy collapse
$Lev(\aleph_0,\lambda)$.  Then, in $V[G]$, there is a tree
isomorphism $h:f(M)\rightarrow f(N)$ as both $f(M)$ and $f(N)$ are
countable and back-and-forth equivalent.  It suffices to prove that
$M$ and $N$ are isomorphic in $V[G]$.

To see this, first note that since `being part of an infinite
branch' is an isomorphism invariant,
the restriction of $h$ to $T_0$ is a tree isomorphism between the $T_0$ of
$M$ and the $T_0$ of $N$.  To finish, we need only show that for every
$\eta\in T_0$ and $j\in\omega$, $M\models P_j(\eta)$ if and only if
$N\models P_j(h(\eta))$.  To see this, let $n=\lg(\eta)$ and $k=\Phi(n,j)$.
Then $M\models P_j(\eta)$ if and only if there is an immediate
successor $\nu$ of $\eta$ that is not part of an infinite branch, but
has an extension $\mu$ of length $n+k$ that is a leaf.  As this condition is
also preserved by $h$, we conclude that $h|_{T_0}$ preserves each of
the $\aleph_0$ colors as well.
\endproof

\begin{Corollary}  There are $\lambda$ pairwise $\bfequiv$-inequivalent
special subtrees of $\lambda^{<\omega}$.
\end{Corollary}

\bp  Let $L=\{R\}$ consist of a single, binary relation, and let
$DG$ be the class of all directed graphs (i.e., $R$-structures) with
universe $\lambda$.  It is well known that there are at least $\lambda$
pairwise $\bfequiv$-inequivalent directed graphs.
But, by composing the maps given in Lemmas
\ref{leftpart} and \ref{rightpart}, we get a $\lambda$-Borel
embedding of $(DG,\bfequiv)$ into $(\hbox{Special subtrees of $\lambda^{<\omega}$},\bfequiv)$ preserving $\bfequiv$ in both directions.
\endproof

\begin{Theorem}  \label{lambdaBorelcomplete}
For any infinite cardinal $\lambda$,
$(\hbox{Subtrees of $\lambda^{<\omega}$},\bfequiv)$ is $\lambda$-Borel complete.
\end{Theorem}

\bp   By Lemma~\ref{leftpart},
it suffices to show $$(\lambda\, CT_\lambda,\bfequiv)\le^B_\lambda
(\hbox{Subtrees of $\lambda^{<\omega}$},\bfequiv)$$
From the Corollary above,
fix a set $\{A_i:i\in\lambda\}$ of pairwise
$\bfequiv$-inequivalent
special subtrees of $\lambda^{<\omega}$.

As notation, let $A_{\<i\>}$ denote the tree $A_i$, and let $A_{\<\>}$
be the two-element tree $\{\<\>,a\}$ satisfying $\<\>\triangleleft a$.
For each $u\subseteq\lambda$, let $T_u=\{\<\>,a\}\cup\{\<i\>:i\in u\}$
and let $S_u=T_u(A_{\<i\>}:i\in u)$.  Note that for each $u\subseteq\lambda$,
$S_u$ has a unique leaf $a$ attached to $\<\>$, and 
the trees $S_u$ and $S_v$ are isomorphic
if and only if $u=v$.

The proof now follows the proof of Lemma~\ref{rightpart}, using the trees
$S_u$ to code the color of a node.  

More formally, let  $T_0:=\lambda^{<\omega}$
and fix an enumeration $\<P_j(x):j\in\lambda\>$ of the unary predicates.
Given any $M\in\lambda\, CT_\lambda$, for each node $\eta\in T_0$,
let $V(\eta):=\{j\in\lambda:M\models V_j(\eta)\}$.
Let $f(M)$ be the tree 
$T_0(S_{V(\eta)}:\eta\in T_0)$.

Note that as each of the $A_i$'s were special, $T_0$ is detectable in
$f(M)$ as being the set of all nodes $\eta$ that are part of an infinite branch
{\bf and} have an immediate successor that is a leaf.
The proof now follows Lemma~\ref{rightpart}.
In particular, given an isomorphism $h:f(M)\rightarrow f(N)$ in $V[G]$,
the restriction of $h$ to $T_0$ is an isomorphism of $M$ and $N$ as
$\kappa\, CT_\lambda$-structures.
\endproof

\section{The Borel completeness of $\aleph_0$-stable, eni-DOP theories}

This section is devoted to the proofs of Theorem~\ref{eniDOPthm} and 
Corollary~\ref{eniDOPcor}.
As the proof of the former is lengthy, the section is split into four subsections.
The first describes two distinct types of eni-DOP witnesses.  The second
shows how one can encode bipartite graphs into models of $T$.  However,
Proposition~\ref{biggy}, which gives a bit of positive information about
the shapes of the bipartite graphs $G$ and $H$ whenever the associated models
$M_G$ and $M_H$ are isomorphic, is rather weak.  Thus, instead of
trying to recover arbitrary bipartite graphs, in the third subsection we
describe how to encode subtrees $\T\subseteq\lambda^{<\omega}$ into bipartite graphs
$G^{[m]}_\T$, where the nodes of $\T$ correspond to complete,
bipartite subgraphs of $G^{[m]}_\T$.
Finally, in the fourth subsection we prove Theorem~\ref{eniDOPthm}, with 
Corollary~\ref{eniDOPcor} following easily from it.

\subsection{Two types of eni-DOP witnesses}  \label{twotypes}

Suppose that $T$ has eni-DOP.  
Call a 5-tuple $(M_0,M_1,M_2,M_3,r)$ an {\em eni-DOP witness} if it satisfies the assumptions of 
Theorem~\ref{equiv}(5).
A {\em finite approximation $\F$} to an eni-DOP witness is a 5-tuple $(a,b,c,d,r_d)$,
where $a,b,c,d$ are finite tuples from $(M_0,M_1,M_2,M_3)$, respectively,
$\tp(b/a)$ and $\tp(c/a)$ are each stationary, regular types,
each of $b,c$ contain $a$ and $\{b,c\}$ are independent over $a$,
$r$ is based and stationary on $d$ with $r_d\in S(d)$ parallel to $r$,
and $\tp(d/bc)\vdash\tp(d/M_1M_2)$.
The last condition, coupled with the fact that $M_0,M_1,M_2$ are each
a-models, yields the following {\em Extendability Condition:}

$$\tp(d/bc)\vdash\tp(d/b^*c^*)$$
for all $a^*\supseteq a$, $b^*\supseteq ba^*$, $c^*\supseteq ca^*$ such that $a^*$
is independent from $bc$ over $a$ and $b^*$ is independent from $c^*$
over $a^*$.  As well, $r_d$ is ENI, $r_d\perp b$, and $r_d\perp c$.

For a fixed choice $\F=(a,b,c,d,r_d)$ of a finite approximation,
the {\em $\F$-candidates over $a$}
consist of all 4-tuples $(b',c',d',r_{d'})$
such that $\tp(a,b,c,d)=\tp(a,b',c',d')$.  There is a natural equivalence
relation $\sim_{\F}$ on the $\F$-candidates over $a$ defined by
$$(b,c,d,r_d)\sim_{\F}(b',c',d',r_{d'})\qquad\hbox{if and only if}
\qquad r_d\not\perp r_{d'}$$

\begin{Lemma} \label{not4}
For any eni-DOP witness $(M_0,M_1,M_2,M_3,r)$,
for any finite approximation $\F$, and for any pair $(b,c,d,r_d)$,
$(b',c',d',r_{d'})$ of equivalent $\F$-candidates over $a$,
every element of the set $\{b,c,b',c'\}$ depends on the other three over
$a$.
\end{Lemma}

\bp
Everything is symmetric, so assume by way of contradiction that
$\fg{b}{a}{cb'c'}$.  First, as $\fg{b'c'}{c}{b}$, the
Extendibility Condition implies that $\tp(d'/b'c')\vdash\tp(d'/b'c'bc)$.
In particular, $\fg{d'}{b'c'}{bc}$, so $\fg{b}{c}{b'c'd'}$ follows by
the symmetry and transitivity of non-forking.

Second, it follows from this and  the Extendibility Condition that
$\tp(d/bc)\vdash\tp(d/bcb'c'd')$, so $\fg{d}{bc}{b'c'd'}$.
Combining these two facts yields
$$\fg{d}{c}{b'c'd'}$$
But then, as $r_d\in S(d)$ is orthogonal to $c$, by e.g., Claim~1.1
of Chapter~X  of \cite{Shc},
$r_d$ would be orthogonal to $b'c'd'$, which contradicts $r_d\not\perp r_{d'}$.
\endproof

It follows from the previous Lemma that there are two types of behavior of
a finite approximation $\F$.  The following definition describes this dichotomy.

\begin{Definition} {\em
Fix an eni-DOP witness $(M_0,M_1,M_2,M_3,r)$.
A finite approximation $\F=(a,b,c,d,r_d)$ of it is {\em flexible\/}
if there is an equivalent $\F$-candidate $(b',c',d',r_{d'})$ over $a$
for which some 3-element subset of $\{b,c,b',c'\}$ is independent over $a$.
We say that the eni-DOP witness $(M_0,M_1,M_2,M_3,r)$ is of
{\em flexible type} if it has a flexible finite approximation.
A witness is {\em inflexible\/} if it is not flexible.
}
\end{Definition}

\begin{Lemma} \label{4.3}
Suppose that $(a,b,c,d,r_d)$ and $(a',b',c',d',r_{d'})$ are each finite approximations of an inflexible
eni-DOP witness satisfying
$\tp(a)=\tp(a')$ and $r_{d}\not\perp r_{d'}$.
Then there is no finite set $A\supseteq aa'$ satisfying $\tp(bc/A)$ does not
fork over $a$, exactly one element from $\{b',c'\}$ is in $A$, and the other element independent from $A$ over $a'$.
\end{Lemma}

\bp  By way of contradiction suppose that $A$ were such a set.
For definiteness, suppose $b'\in A$ and $\fg{c'}{a'}{A}$.
Let $\F$ denote the finite approximation exemplified by
$(A,bA,cA,dA,r_{dA})$.
Fix an automorphism $\sigma\in Aut(\C)$ fixing $Ac'$ pointwise
such that $\fg{bcd}{Ac'}{\sigma(b)\sigma(c)\sigma(d)}$.
Then $(\sigma(b)A,\sigma(c)A,\sigma(d)A,r_{\sigma(d)A})$
is an $\F$-candidate over $A$.
Moreover,
since $r_{dA}\not\perp r_{d}\not\perp r_{d'}\not\perp r_{\sigma(d)A}$
the transitivity of non-orthogonality of regular types imply that
 it is equivalent to $(bA,cA,dA,r_{dA})$.
We will obtain a contradiction to the inflexibility of the eni-DOP witness
by exhibiting a 3-element subset of $\{b,c,\sigma(b),\sigma(c)\}$
that is independent over $A$.

To see this, first note that since $b$ and $c$ are independent over $A$
and $\tp(c'/A)$ has weight 1, $c'$ cannot fork with both $b$ and $c$ over $A$.
For definiteness, suppose that
$b$ and $c'$ are independent over $A$.
It follows that $\sigma(b)$ is also independent from $c'$ over $A$.
These facts, together with the
independence of $b$ and $\sigma(b)$ over $Ac'$, imply that the
three element set $\{b,\sigma(b),c'\}$ is independent over $A$.

We next claim that $\tp(bc/Ac')$ forks over $A$.  If this were not the case,
recalling that $b'\in A$,
we would have $\fg{bc}{aa'}{b'c'}$.
Then, by two applications of the Extendibility Condition,
we would have $\fg{bcd}{aa'}{b'c'd'}$, which would
contradict $r_d\not\perp r_{d'}$.

But now, the results in the previous two paragraphs, together with the
fact that $\tp(c/Ab)$ has weight 1, imply that the set $\{b,\sigma(b),c\}$
is independent over $A$, contradicting the inflexibility of the
eni-DOP witness.
\endproof

\subsection{Coding bipartite graphs into models}  \label{4.2}

In this subsection, we take a particular eni-DOP witness and show
how we can embed an arbitrary bipartite graph $G$ into a model $M_G$.
This mapping will be Borel, and isomorphic graphs will give rise to
isomorphic models, but the converse is less clear.  Proposition~\ref{biggy}
demonstrates that the graphs $G$ and $H$ must be similar in some weak sense
whenever $M_G$ and $M_H$ are isomorphic.

Fix an eni-DOP witness $(M_0,M_1,M_2,M_3,r)$
and a finite approximation $\F=(a,b,c,d,r_d)$ of it,
choosing $\F$ to be flexible if the witness is.
As notation, let $p=\tp(b/a)$ and $q=\tp(c/a)$.

We begin by describing
 how to code arbitrary bipartite graphs into models of $T$.
Given a bipartite graph $G=(L_G,R_G,E_G)$,
choose sets $\B_G:=\{b_g:g\in L_G\}$ and $\CC_G:=\{c_h:h\in R_G\}$
such that $\B_G\cup\CC_G$ is independent over $a$,
$\tp(b_g/a)=p$ for each $b_g\in\B_G$, and $\tp(c_h/a)=q$ for each
$c_h\in\CC_G$.  As well, for each $(g,h)\in L_G\times R_G$, choose an element
$d_{g,h}$ such that $\tp(d_{g,h}b_gc_h/a)=\tp(dbc/a)$
and let $r_{g,h}\in S(d_{g,h})$ be conjugate to $r_d$.
Note that $r_{g,h}\perp r_{g',h'}$ unless $(g,h)=(g',h')$.
Let $\D_G=\{d_{g,h}:(g,h)\in E_G\}$ and let
$\R_G=\{r_{g,h}:(g,h)\in E_G\}$.

Inductively construct models $M^n_G$ of $T$ as follows:
$M^0_G$ is any prime model over $\B_G\cup\CC_G\cup\D_G$.
Given $M^n_G$, let $\Pcal_n=\{p\in S(M^n_G):p\perp \R_G\}$.
By the $\aleph_0$-stability of $T$, $\Pcal_n$ is countable.
Let $\E_n=\{e_s:s\in \Pcal_n\}$ be independent over $M^n_G$
with each $e_s$ realizing $s$, and let $M^{n+1}_G$
be prime over $M^n_G\cup\E_n$.
Finally, let $M_G=\bigcup_{n\in\omega} M^n_G$.

It is easily verified that  if $G$ has universe $\lambda$,
then the mapping $G\mapsto M_G$ is $\lambda$-Borel.
Moreover, it is easy to see that for regular types $r\in S(M_G)$,
$$r \ \hbox{has finite dimension in $M_G$ if and only if $r\not\perp r_{g,h}$
for some $(g,h)\in E_G$}$$

Suppose that $f:M_G\rightarrow M_H$ were an isomorphism.
Then $f$ maps the regular types in $S(M_G)$ of finite dimension
onto the regular types in $S(M_H)$ of finite dimension.
Thus, by construction of $M_G$ and $M_H$, this correspondence
yields a bijection
$$\pi_f:E_G\rightarrow E_H$$
Unfortunately, this identification need not extend to a bipartite graph
isomorphism between $G$ and $H$.  Specifically, there might be
edges $e_1,e_2\in E_G$ that share a vertex of $G$, while
the corresponding edges $\pi_f(e_1),\pi_f(e_2)\in E_H$ do not
have a common vertex.   The bulk of our argument will be to show
that images of sufficiently large, complete
bipartite subgraphs of $G$ cannot be too wild.

To make this precise, for $X\subseteq E_G$, let $v_G(X)$ denote
the smallest subset of the vertices of $G$ with $X\subseteq E_{v_G(X)}$.
For $\ell$ very large, call a graph $G$ {\em almost $\ell$-complete bipartite}
if it is $m_1\times m_2$ bipartite with $0.99\ell\le m_i\le \ell$ for $i=1,2$
and each vertex has valence at least $0.9\ell$.

The proof of the following Proposition is substantial, and occupies the remainder
of this subsection.

\begin{Proposition} \label{biggy}
For any bipartite graphs $G$ and $H$
and for any isomorphism $f:M_G\rightarrow M_H$,
there is a number $\ell^*$, depending only on $f$, such that for all
$\ell\ge\ell^*$, if $G_0\subseteq G$ is any complete $\ell\times\ell$ bipartite
subgraph, then $v_H(\pi_f(E_{G_0}))$ contains an almost $\ell$-complete bipartite subgraph.
\end{Proposition}

\bp
Fix bipartite graphs $G$, $H$, and an isomorphism
$f:M_G\rightarrow M_H$.  As notation, let $a'=f^{-1}(a)$,
let $\B_H'=\{f^{-1}(b):b\in \B_H\}$, and let $\CC_H'=\{f^{-1}(c):c\in\CC_H\}$.
Let $X\subseteq \B_G\cup\CC_G$ be minimal such that
$\tp(a'/a\B_G\cup\CC_G)$ does not fork over $Xa$,
and let $X'\subseteq\B_H'\cup\CC_H'$ be minimal such that
$\tp(a/a'\B_H'\CC_H')$ does not fork over $X'a'$.
Note that $|X|\le \wt(a'/a)$ and $|X'|\le\wt(a/a')$.

Let $\Lambda^*$ be the set of non-orthogonality classes of
regular types in $S(M_G)$ of finite dimension in $M_G$.
For each $S\in\Lambda^*$ let $(b_s,c_s)$ be the unique element of
$\B_G\times\CC_G$ such that
there is a candidate $(a,b_s,c_s,d,r_d)$ over $a$
with $r_d\in S$ and let $(b_s',c_s')$ be the unique element of
$\B_H'\times\CC_H'$ such that there is a
 candidate $(a',b_s',c_s',d',r_{d'})$ over $a'$ with
$r_{d'}\in S$.

For $\Lambda$ a finite subset of $\Lambda^*$, let
$B(\Lambda)=\{b_s:S\in\Lambda\}$, $C(\Lambda)=\{c_s:S\in \Lambda\}$,
and $v(\Lambda)=B(\Lambda)\cup C(\Lambda)$.
Dually, define $B'(\Lambda)$, $C'(\Lambda)$, and $v'(\Lambda)$
using $(b'_s,c'_s)$ in place of $(b_s,c_s)$.

The proof splits into two cases
depending on whether our eni-DOP witness is flexible or inflexible.

\medskip\par
\noindent{\bf Case 1:}  The eni-DOP witness is inflexible.
\medskip

This case will be substantially easier than the other, and in fact, we
prove that there is a number $e$ such that for all sufficiently large $\ell$,
the image of any $\ell\times\ell$ bipartite graph contains an $(\ell-e)\times(\ell-e)$
complete, bipartite subgraph.  The simplicity of this case is
primarily due to the following claim.

\medskip\par
\noindent{\bf Claim 1.}  For any finite $\Lambda\subseteq\Lambda^*$
such that $v(\Lambda)$ is disjoint from $X$ and $v'(\Lambda)$
is disjoint from $X'$, we have $|v(\Lambda)|=|v'(\Lambda)|$.

\bp To see this, we again split into cases.  First, if $p\perp q$,
then we handle the two `halves' separately.
Note that for each $S\in \Lambda$, $\tp(b_sc_s/aa')$ does not fork
over $a$, $\tp(b'_s,c'_s/aa')$ does not fork over $a'$, and by Lemma~\ref{4.3},
each element of $\{b_s,c_s\}$ forks with $b'_sc'_s$ over
$aa'$.  Since $p\perp q$, this implies $\{b_s,b'_s\}$ fork over $aa'$.
It follows that, working over $aa'$,
$$Cl_p(B(\Lambda))=Cl_p(B'(\Lambda))$$ hence
$|B(\Lambda)|=|B'(\Lambda)|$.  It follows by a symmetric argument that
$Cl_q(C(\Lambda))=Cl_q(C'(\Lambda))$, so $|C(\Lambda)|=|C'(\Lambda)|$.
It follows immediately that $|v(\Lambda)|=|v'(\Lambda)|$.

On the other hand, if $p\not\perp q$, then $Cl_p$ is a closure relation on
$p^*(\C)\cup q^*(\C)$, where $p^*$ (resp.\ $q^*$) is the non-forking
extension of $p$ (resp.\ $q$) to $aa'$.
Furthermore, for each $S\in\Lambda$ we have
$Cl_p(b_sc_s)=Cl_p(b'_s,c'_s)$.  It follows that
$Cl_p(v(\Lambda))=Cl_p(v'(\Lambda))$.  As each set is independent over
$aa'$, we conclude that $|v(\Lambda)|=|v'(\Lambda)|$.

Let $w=\wt(a'/a)$ and $e=w+\wt(a/a')^2$.
Suppose that $G_0\subseteq G$ is an $\ell\times\ell$ complete,
bipartite subgraph.
Since $|X|\le w$, there is an $(\ell-w)\times(\ell-w)$  complete subgraph
$G_0^*\subseteq G_0$ such that $E_{G_0^*}$ is disjoint from $X$.
By our choice of $e$ there is an $(\ell-e)\times(\ell-e)$ complete
subgraph $G_1\subseteq G_0^*$ such that $\pi_f(b,c)$ is not contained
in $X'$ for all pairs $(b,c)\in E_{G_1}$.  But then, by
Lemma~\ref{4.3}, we have $\pi_f(b,c)$ is disjoint from $X'$ for
all $(b,c)\in E_{G_1}$.

Now, $G_1$ is an $(\ell-e)\times(\ell-e)$ complete, bipartite subgraph of
$G$.  In particular, $G_1$ has $2(\ell-e)$ vertices and $(\ell-e)^2$
edges.  Let $H_1$ be the subgraph of $H$ whose edges are $E_{H_1}:=\pi_f(E_{G_1})$
and whose vertices are $v(H_1):=v_H(E_{H_1})$.
Then $|E_{H_1}|=(\ell-e)^2$ since $\pi_f$ is a bijection and
$$|v(H_1)|=|v_H(E_{H_1})|=|v_G(E_{G_1})|=2(\ell-e)$$
by Claim 1.  By a classical optimal packing result, this is only possible
when $H_1$ is itself a complete, $(\ell-e)\times(\ell-e)$ bipartite subgraph of
$H$.

\medskip\par
\noindent{\bf Case 2:}  The eni-DOP witness is flexible.
\medskip

As we insisted that our finite approximation be flexible, it follows  from
Lemma~\ref{not4} that $p\not\perp q$, so $p$-closure is a
dependence relation on $p(\C)\cup q(\C)$.

As well, for any candidate $(b,c,d,r_d)$ over $a$ and for any finite
$A\supseteq a$, there is an equivalent candidate
$(b',c',d',r_{d'})$ over $a$ such that
$w_p(b'c'/A)=1$.

\begin{Definition} {\em  For any finite subgraph $G_0\subseteq G$, let
$\Lambda(G_0)$ be the set of non-orthogonality classes
$\{[r_{d_{g,h}}]: (g,h)\in E_{G_0}\}.$
Note that $|\Lambda(G_0)|=|E_{G_0}|$ by the pairwise orthogonality of the types
$r_{d_{g,h}}$.

A {\em manifestation $\M=\M(\Lambda,a)$ over $a$} is a set of
candidates $\{(b_s,c_s,d_s,r_{d_s}):S\in\Lambda\}$  over $a$ with
$r_{d_s}\in S$ for each $S\in\Lambda$.  Associated to any manifestation $\M$ is a bipartite
graph $G(\M)$ with `Left Nodes' $L(\M)=\{b_s:s\in\Lambda\}$,
`Right Nodes' $R(\M)=\{c_s:s\in\Lambda\}$, vertices
$v(\M)=L(\M)\cup R(\M)$, and edges
$E(\M)=\{(b_s,c_s):s\in\Lambda\}$.

If $G_0$ is a subgraph of $G$, then the {\em canonical manifestation
of $\Lambda(G_0)$ over $a$ inside $M_G$\/} is the set
$$\{(b_g,c_h,d_{g,h},r_{g,h}):(g,h)\in E_{G_0}\}$$.

A set $A$ {\em represents $\Lambda$ over $a$} if
$a\subseteq A$ and $v(\M)\subseteq A$ for some manifestation
$\M$ of $\Lambda$ over $a$.  A manifestation $\M'$ is
{\em $A$-free} if
$w_p(b_s',c_s'/A)=1$ for each $S\in\Lambda$
and $\{(b_s',c_s'):S\in\Lambda\}$ are independent over $A$.

}
\end{Definition}

Now, working in the monster model $\C$, we define a measure
of the complexity of $\Lambda$ over $a$.
First, note that for any candidate $(b,c,d,r_d)$ over $a$,
there is an equivalent candidate $(b','c',d',r_{d'})$ over $a$
with $w_p(b'c'/abc)=1$.  By choosing $b'c'$ to be independent
over $abc$ from
any given $A\supseteq abc$ we can insist that $w_p(b'c'/A)=1$.
It follows that $A$-free manifestations of $\Lambda$ exist over
any set $A$ representing a finite $\Lambda$.
Thus, the following definition makes sense.

\begin{Definition} {\em  The {\em maximal weight,\/}
 $mw(\Lambda,a)$, is the largest integer $m$
such that for all finite $A$ representing $\Lambda$ over $a$,
there is an $A$-free manifestation $\M'(\Lambda,a)$ over $a$ with
$|v(\M')|=m+\Lambda$.
}
\end{Definition}

\begin{Lemma}  \label{connected}
Suppose that $G$ is a bipartite graph, $G_0\subseteq G$ is a
connected subgraph of $G$, let $\M(\Lambda(G_0),a)$ be the canonical
manifestation of $\Lambda(G_0)$ inside $M_G$, and let
$\M'(\Lambda,a)$ be any other manifestation of $\Lambda(G_0)$.
Then $$Cl_p(v(\M')\cup\{v\})=Cl_p(v(\M')\cup v(G_0))$$
for any $v\in v(G_0)$.
\end{Lemma}

\bp  Arguing by symmetry and induction,
it suffices to show that for all nonempty
$B\subseteq v(G_0)$ and every $c\in v(G_0)\setminus B$ such that
$(b,c)\in E_{G_0}$ for some $b\in B$ we have
that $c\in Cl_p(v(\M')\cup B)$.  But this is immediate, since
$Cl_p(\{b',c',b,c\})=Cl_p(\{b',c',b\})$ for all equivalent
candidates $(b,c,d,r_d)$ and $(b',c',d',r_{d'})$ over $a$.
\endproof

\begin{Lemma} \label{leftright}
$k(G_0)\le mw(\Lambda(G_0),a)\le |v(G_0)|$
for any bipartite graph $G$ and any finite $G_0\subseteq G$.
\end{Lemma}

\bp  The upper bound is very soft.  Let $A\supseteq a\cup v(G_0)$
be arbitrary and let $\M'$ be any other manifestation of
$\Lambda(G_0)$ over $a$.
Then $$w_p(v(\M')/a)\le w_p(v(\M')v(G_0)/a)=
w_p(v(\M')/av(G_0)) + w_p(v(G_0)/a)$$
Since $w_p(b'_sc'_s/ab_sc_s)\le 1$ for each $S\in \Lambda(G_0)$,
we have $w_p(v(\M')/av(G_0))\le |\Lambda(G_0)|$.
Also, by the independence of the nodes in $M_G$, $w_p(v(G_0)/a)=|v(G_0)|$.
The upper bound on $mw(\Lambda(G_0),a)$ follows immediately.

For the lower bound, again choose any $A\supseteq av(G_0)$
and let $C\subseteq v(G_0)$ consist of one vertex from every
connected component of $G_0$.  Clearly, $A$ represents $\Lambda(G_0)$ and
$|C|=CC(G_0)$.
Let $\M'$ be any $A$-free manifestation of $\Lambda(G_0)$ over $a$.
Then
$$w_p(v(\M')/a)\ge w_p(v(\M')C/a)-CC(G_0)=
w_p(v(\M')v(G_0)/a)-CC(G_0)$$
with the second equality coming from Lemma~\ref{connected}.
As before, for each $S\in\Lambda(G_0)$,
$w_p(b_s'c_s'/ab_sc_s)\le 1$ so $w_p(v(\M')/av(G_0))\le |\Lambda(G_0)|$.
On the other hand, the $A$-freeness of $\M'$ implies that
$w_p(v(\M')/A)=|\Lambda(G_0)|$, hence
$w_p(v(\M')/av(G_0))=|\Lambda(G_0)|$.
Thus, $$w_p(v(\M')v(G_0)/a)=w_p(v(\M')/v(G_0)a)+w_p(v(G_0)/a)
=|\Lambda(G_0)|+|v(G_0)|$$
from which the lower bound follows as well.
\endproof

Now, returning to our isomorphism $f:M_G\rightarrow M_H$,
suppose that $G_0$ is any finite subgraph of $G$ that is disjoint from $X$,
i.e., so that $\tp(G_0/aa')$ does not fork over $a$.  We then claim:

\medskip\par
\noindent{\bf Claim 2:}  $mw(\Lambda(G_0),a')\le |v(G_0)| + wt(a/a')$

\bp  Choose any finite $A$ containing $\{aa'\}\cup v(G_0) \cup v_H(\pi_f(E_{G_0}))$.
So $A$ represents $\Lambda(G_0)$ over $a'$.
Let $\M'$ be any $A$-free manifestation of $\Lambda(G_0)$ over $a'$.
Now
$$w_p(v(\M')/a')\le w_p(v(\M')aG_0/a')=w_p(v(\M')/aa'G_0)+w_p(aG_0/a')$$
But, as before $w_p(b'_sc'_s/aa'b_sc_s)\le 1$, so
$w_p(v(\M')/aa'G_0)\le |\Lambda(G_0)|$.  Also,
$$w_p(aG_0/a')=w_p(G_0/aa')+wt(a/a')=|v(G_0)|+w_p(a/a')$$
and the Claim follows.
\endproof

Finally, choose a complete, bipartite subgraph $G_0\subseteq G$,
where $\ell$ is sufficiently large with respect to $W=wt(a/a')$.
Let $H_0$ be the bipartite graph with vertices $v_H(\pi_f(E_{G_0}))$
and edges $\pi_f(E_{G_0})$ and let $H_0^*$ be the subgraph of $H$ with
the same vertex set as $H_0$.  Note that $E_{H_0}\subseteq E_{H_0^*}$,
but that equality need not hold.

As $G_0$ is $\ell\times\ell$ complete bipartite, $|v(G_0)|=2\ell$ and
$|\Lambda(G_0)|=\ell^2$.  It follows immediately that $|E_{H_0}|=\ell^2$
and it follows from Claim~2 and Lemma~\ref{leftright} that
$$k(H_0)\le mw(\Gamma(G_0),a')\le 2\ell +W$$
where $W=wt(a/a')$.
So, by Corollary~\ref{finitecomb}
of the Appendix, $H_0$ contains an almost $\ell$-complete
bipartite subgraph $H_1$.  But then, $H_1^*$, which is the subgraph of $H$ with the same
vertex set as $H_1$, is almost $\ell$-complete as well.
\endproof

\subsection{Coding trees by complete, bipartite subgraphs}  \label{sec4.3}

As Proposition~\ref{biggy} is rather weak, we give up on coding arbitrary bipartite graphs
into models of $T$.  Rather, we seek to code subtrees of $\lambda^{<\omega}$ into bipartite
graphs that have large, complete subgraphs.

Fix a sufficiently large integer $m$ and a
tree $\T\subseteq \lambda^{<\omega}$.
We will construct a bipartite graph $G^{[m]}_\T$,
whose  $7m\times 7m$ complete bipartite subgraphs $B^m_\T(\eta)$
 code nodes  $\eta\in \T$.
Moreover, additional information about
the level of $\eta$ and its set of immediate successors will be coded
by the size of the intersection of $B^m_\T(\eta)$ and $B^m_\T(\nu)$ for other
$\nu\in T$.

More precisely, fix a tree $(\T,\trianglelefteq)$ and a large integer $m$.
We first define a bipartite graph $preG^{[m]}_\T$ to have universe
$\T\times m\times 14$ with the edge relation
$$\{((\eta,i_1,n_1),(\eta,i_2,n_2)):\eta\in \T, i_i,i_2\in m,
 \hbox{$n_1+n_2$ is odd}\}$$
So the `left hand side' of $pre G^{[m]}_\T$ is
$L=\T\times m\times\{n\in 14:n$ odd$\}$,
the `right hand side' is
$R=\T\times m\times\{n\in 14:n$ even$\}$,
thereby associating a $7m\times 7m$ complete, bipartite graph
to each node $\eta\in \T$.

Next, define a binary relation $E_0$ on $preG^{[m]}_\T$ by
$(\eta_1,i_1,n_1) E_0 (\eta_2,i_2,n_2)$ if and only if
\begin{itemize}
\item $\eta_2$ is an immediate successor of $\eta_1$, $i_1=i_2$, $n_1=n_2$
and
\item either
$\lg(\eta_1)=0$ and $n_1\in\{0,1\}$ or
$\lg(\eta_1)>0$ and $n_1\in\{10,11,12,13\}$.
\end{itemize}

Let $E$ be the smallest equivalence relation containing $E_0$,
i.e., the  reflexive, symmetric and transitive closure of $E_0$.

Let $G^{[m]}_\T:=preG^{[m]}_\T/E$
and, for each $\eta\in\T$, let $B^m_\T(\eta)=\{g\in G^{[m]}_\T:
(\eta,i,n)\in g$ for some $i<m,n<14\}$.

As notation, for each $\eta\in\T$, let $B^m_\T(\eta)=\{g\in G^{[m]}_\T:
(\eta,i,n)\in g$ for some $i<m,n<14\}$, let
$\S^m_\T=\{B^m_\T(\eta):\eta\in\T\}$,
and let $g_\T:\T\rightarrow\S^m_\T$ be the bijection
$\eta\mapsto B^m_\T(\eta)$.  For all of these, when $\T$ and $m$ are
clear, we delete reference to them.  Finally, call an element $g\in G^{[m]}_\T$
a {\em singleton} if $g=\{(\eta,i,n)\}$ for a single element
$(\eta,i,n)\in preG^{[m]}_\T$.
All of the following Facts are immediate:

\begin{Fact} \label{gg}
\begin{enumerate}
\item Every $B(\eta)$ is a $7m\times 7m$ complete, bipartite graph;
\item  If $g\in B(\eta)$ is a singleton and $E(g,h)$, then
$h\in B(\eta)$;
\item  For all $\eta\in\T$, $i<m$, $(\eta,i,n)$ is a singleton for all
$2\le n\le 9$.
\item \label{important}
 If $lg(\nu)<lg(\eta)$, then $B(\nu)\cap B(\eta)=\emptyset$ if and only if
$\nu=\eta^-$.  Moreover, a nonempty intersection is a complete $m\times m$ bipartite graph if $\eta^-=\<\>_\T$ and the intersection is
$2m\times 2m$ complete, bipartite if $\eta^-\neq\<\>_\T$.
\end{enumerate}
\end{Fact}

\begin{Lemma}  \label{onlycomplete}
$\S=\{all$ $7m\times 7m$ complete, bipartite subgraphs of $G^{[m]}_\T\}$.
\end{Lemma}

\bp  That each $B(\eta)\in\S$ is a $7m\times 7m$ complete, bipartite
subgraph is clear.  Conversely, fix a  $7m\times 7m$ complete, bipartite
subgraph of $G^{[m]}_\T$.  First, suppose that $X$ contains a singleton
$a$.  Without loss, assume $a\in X\cap B(\eta)\cap L$.
Then $E_X(a)=\{b\in X:E(a,b)\}$ has cardinality $7m$ and is contained
in $B(\eta)\cap R$, hence $E_X(a)=B(\eta)\cap R$.  But then, $X\cap R$
contains a singleton as well, so arguing similarly, $B(\eta)\cap L=X\cap L$,
so $X=B(\eta)$.
It remains to show that $X$ contains a singleton.  Choose $k$ maximal
such that there is $\eta\in\T$, $lg(\eta)=n$, and $X\cap B(\eta)\neq\emptyset$.
Let $\nu=\eta^-$.  If $X$ does not contain a singleton, then the maximality
of $k$ implies that $X\cap(B(\eta')\setminus B(\nu))=\emptyset$ for
all $\eta'\in Succ(\nu)$.  Choose any $a\in X\cap B(\eta)\cap L$.
Then $a\in B(\nu)$ and moreover, $E_X(a)\subseteq B(\nu)\cap R$.
By counting, $E_X(a)=B(\nu)\cap R$, so $X$ contains a singleton,
completing the proof of the Claim.
\endproof

For clarity, let $L_0=\{R_1,R_2\}$ denote the language consisting of two
binary relation symbols.  Form an $L_0$ structure
$(\S^m_\T,R_1,R_2)$ by declaring that $R_1(X,Y)$ holds if and only
if $X\cap Y$ is an $m\times m$ complete, bipartite graph and
$R_2(X,Y)$ holds if and only if $X\cap Y$ is a $2m\times 2m$
complete, bipartite graph.

\begin{Lemma} \label{suff} For any sufficiently large $m$ and trees
$(\T,\trianglelefteq), (\T',\trianglelefteq)$, if
there is an $L_0$-isomorphism $\Phi:(\S^m_\T,R_1,R_2)\rightarrow
(\S^m_{\T'},R_1,R_2)$ of the associated $L_0$-structures,
then the composition $h:(\T,\trianglelefteq)\rightarrow (\T',\trianglelefteq)$
given by $h=g_{\T'}^{-1}\circ \Phi\circ g_\T$ is a tree isomorphism.
\end{Lemma}

\bp  For each $n\in\omega$, let $\T_n=\{\eta\in\T:lg(\eta)<n\}$
and define $\T'_n$ analogously.  Using Fact~\ref{gg}(\ref{important}),
one proves by induction on $n$ that
$h|_{\T_n}:(\T_n,\trianglelefteq)\rightarrow (\T'_n,\trianglelefteq)$
is a tree isomorphism.  This suffices to prove the Lemma.
\endproof

\subsection{$\aleph_0$-stable, eni-DOP theories are $\lambda$-Borel complete}

\begin{Theorem}  \label{eniDOPthm}
If $T$ is $\aleph_0$-stable with eni-DOP, then 
for any infinite cardinal $\lambda$, there is a $\lambda$-Borel embedding
$\T\mapsto M(\T)$ from subtrees of $\lambda^{<\omega}$ to $Mod_\lambda(T)$
satisfying 
$$(\T_1,\trianglelefteq)\cong (\T_2,\trianglelefteq)\qquad\hbox{if and only if}\qquad M(\T_1)\cong M(\T_2)$$
\end{Theorem}

\bp  Fix any infinite cardinal $\lambda$.  
As in Subsection~\ref{4.2}, 
fix an eni-DOP witness $(M_0,M_1,M_2,M_3,r)$
and a finite approximation $\F=(a,b,c,d,r_d)$ of it,
choosing $\F$ to be flexible if the witness is.
As notation, let $p=\tp(b/a)$ and $q=\tp(c/a)$.
As well, for the whole of the proof, fix a recursive, fast growing sequence,
$\<m_i:i\in\omega\>$ of integers, e.g., $m_0=10$ and $m_{i+1}=m_i!!$

Given a subtree $\T\subseteq \lambda^{<\omega}$, let $G^*_\T$ be the bipartite graph
which is the disjoint union $\bigcup_{i\in\omega} G^{[m_i]}_\T$, where the graphs
$G^{[m_i]}_\T$ are constructed as in Subsection~\ref{sec4.3}.
Next, construct a model $M(\T):=M_{G^*_\T}$ from the bipartite graph $G^*_\T$ as in 
Subsection~\ref{4.2}.
Clearly, after some reasonable coding, we may assume that $M(\T)$ has universe
$\lambda$.  It is routine to verify that
both of the maps $\T\mapsto G^*_\T$ and $G^*_\T\mapsto
M_{G^*_\T}$ (and hence their composition) are $\lambda$-Borel.

By looking at the constructions in Subsections~\ref{4.2} and \ref{sec4.3}, it is
easily checked isomorphic trees $\T\cong\T'$ give rise to isomorphic models
$M(\T)\cong M(\T')$.

To establish the converse, suppose that $\T,\T'$ are subtrees
such that $M(\T)\cong M(\T')$.  
Fix an isomorphism $f:M_{G^*_{\T}}\rightarrow
M_{G^*_{\T'}}$ and choose $i$ so that $m_i>>\ell^*$, where $\ell^*$
is the constant in the statement of Proposition~\ref{biggy}.

For each $\eta\in\T$, by Fact~\ref{gg}(1), $B^m_\T(\eta)$
is a $7m_i\times 7m_i$ complete, bipartite subgraph of $G^{[m_i]}_\T$.
Let $E(\eta)$ denote the edge set of $B^m_\T(\eta)$.
Now $\pi_f(E(\eta))$ is a set of $(7m_i)^2$ edges in $G^*_{\T'}$.
Let $v_{\T'}(\eta)$ be the smallest set of vertices in $G^*_{\T'}$
whose edge set contains $\pi_f(E(\eta))$.

By Proposition~\ref{biggy}, the graph
 $J(\eta):=(v_{T'}(\eta),\pi_f(E(\eta)))$ has an almost $7m_i$-complete
bipartite subgraph $K(\eta)$.  Let $K^*(\eta)$ be the subgraph of
$G^*_{\T'}$ whose vertex set is the same as $K(\eta)$.  Note
that the edge set of $K^*(\eta)$ contains the edge set of $K(\eta)$,
so $K^*(\eta)$ is almost $7m_i$-complete as well.

As $K^*(\eta)$ is a connected subgraph of $G^*_{\T'}$,
$K^*(\eta)\subseteq G^{[m_j]}_{\T'}$ for some $j$.
As the valence of each vertex of $K^*(\eta)$ is $\sim 7m_i$
and $m_i>>m_k$ for all $k<i$, we must have $j\ge i$.

\medskip\par\noindent{\bf Claim:}  $j=i$.

\medskip\par
\bp  Choose $\nu\in\T'$ such that $K^*(\eta)$ and $B^{m_j}(\nu)$
share a connected subgraph $D$  with $e(D)>>N_f$.
Arguing as above, there is an almost $7m_j$ complete, bipartite
subgraph $H^*(\nu)$ of $G^*_\T$ whose edge set (almost) contains
$\pi_f^{-1}(E(\nu))$, where $E(\nu)$ is the edge set of $B^{m_j}_{\T'}(\nu)$.
As before, $H^*(\nu)\subseteq G^{[m_k]}_{\T}$ for some $k$,
and as the valence of every vertex is large, $k\ge j$.
However, almost all of the edges of $D$ correspond to edges
of $H^*(\nu)$.  In particular, $H^*(\nu)$ contains edges from
$B^{m_i}_\T(\eta)$.  But, as $H^*(\eta)$ is connected, this implies
$H^*(\eta)\subseteq G^{[m_i]}_\T$.  Thus $k=j=i$.
\endproof

Thus, we have shown that for each $\eta\in\T$, $K^*(\eta)$
is an almost $7m_i$ complete bipartite subgraph of $G^{[m_i]}_{\T'}$.
It follows as in the proof of Lemma~\ref{onlycomplete} that for each
$\eta\in\T'$, there is a unique $\nu\in\T'$ such that
the subgraphs $K^*(\eta)$ and $B^{m_i}_{\T'}(\nu)$ have large
intersection in $G^{[m_i]}_{\T'}$.
Define $$\Phi:\S^{m_i}_\T\rightarrow \S^{m_i}_{\T'}$$
by $\Phi(B^{m_i}(\eta))=B^{m_i}_{\T'}(\nu)$ for this unique $\nu$.
As the argument given above is reversible, $\Phi$ is a bijection.
Furthermore, if $D\subseteq B^{m_i}(\eta)$ is either an $m\times m$
or a $2m\times 2m$ complete, bipartite subgraph, then applying
Proposition~\ref{biggy} to $D$ yields a connected graph
$K^*(D)$  whose number of edges satisfies
$$m_i^2-N_f\le e(K(D))\le m_i^2$$
By taking $D=B^{m_i}(\eta_1)\cup B^{m_i}(\eta_2)$
for various $\eta_1,\eta_2\in\T$, it follows that $\Phi$
is an $L_0$-isomorphism.  Thus, by Lemma~\ref{suff},
$(\T,\trianglelefteq)\cong(\T',\trianglelefteq)$ as required.
\endproof

\begin{Corollary} \label{eniDOPcor}
If $T$ is $\aleph_0$-stable with eni-DOP, then 
$T$ is Borel complete.  Moreover, for every infinite cardinal $\lambda$,
$T$ is $\lambda$-Borel complete for $\bfequiv$.
\end{Corollary}

\bp  For both statements, by Theorem~\ref{lambdaBorelcomplete}, it suffices to 
show that
$$(\hbox{Subtrees of $\lambda^{<\omega}$},\bfequiv)\ \le^B_\lambda  \
(Mod_\lambda(T),\bfequiv)$$
for every $\lambda\ge\aleph_0$.
So fix an infinite cardinal $\lambda$.  The map $\T\mapsto M(T)$ given
in Theorem~\ref{eniDOPthm} is $\lambda$-Borel.  Choose any generic filter $G$ for the
Levy collapsing poset $Lev(\lambda,\aleph_0)$.
By Fact~\ref{Levy}, for any subtrees $\T_1,\T_2\subseteq\lambda^{<\omega}$ in $V$,
$\T_1\bfequiv\T_2$ in $V$ if and only if $\T_1\cong\T_2$ in $V[G]$.
As well, $M(\T_1)\bfequiv M(\T_2)$ in $V$ if and only if $M(\T_1)\cong M(\T_2)$ in
$V[G]$.  Thus, since the mapping $\T\mapsto M(\T)$ is visibly absolute between $V$ and $V[G]$,
the result follows immediately from Theorem~\ref{eniDOPthm}.
\endproof

\section{eni-NDOP and decompositions of models}

In this section, we assume throughout that $T$ is $\aleph_0$-stable
with eni-NDOP. [In fact, the first few Lemmas require only $\aleph_0$-stability.]
 We discuss three species of decompositions (regular,
eni, and eni-active) of an arbitrary model $M$ and prove a
theorem about each one.   Theorem~\ref{regulartheorem} asserts
that in a regular decomposition $\d=\<M_\eta,a_\eta:\eta\in I\>$ of $M$,
then $M$ is atomic over $\bigcup_{\eta\in I} M_\eta$.  This theorem
plays a key role in Corollary~\ref{notop}.

Next, we discuss eni-active decompositions of a model $M$ and prove
that for any $N\preceq M$ that contains $\bigcup_{\eta\in I} M_\eta$,
then $N$ is an $L_{\infty,\aleph_0}$-elementary substructure of $M$.
In particular, Corollary~\ref{quotelater} states that an eni-active
decomposition determines a model up to $L_{\infty,\aleph_0}$-equivalence.
This is extremely important when we compute 
$I_{\infty,\aleph_0}(T,\kappa)$ in Section~\ref{last}.

Finally, we prove Theorem~\ref{enitheorem}, which states that
a model $M$ is atomic over $\bigcup_{\eta\in I} M_\eta$ for any
eni decomposition of $M$ {\em provided that each of the models is
maximal atomic} (see Definition~\ref{maxatomic}).  While the result
sounds strong, it is of little use to us, as one has little control about what
the maximal atomic submodels of an arbitrary model look like.
This theorem was also proved by Koerwien~\cite{K}, but is included here
to contrast with Theorems~\ref{regulartheorem} and \ref{eniactivetheorem}.

\begin{Definition} {\em  An {\em independent tree of models}
$\{M_\eta:\eta\in I\}$ satisfies
\begin{itemize}
\item $I$ is a subtree of $Ord^{<\omega}$;
\item $\eta\trianglelefteq\nu$ implies $M_\eta\preceq M_\nu$; and
\item  For each $\eta\in I$ and
$\nu\in Succ_I(\eta)$, $\fg {\bigcup_{\nu\trianglelefteq\gamma} M_\gamma}
{M_\eta} {\bigcup_{\nu\not\trianglelefteq\delta}M_\delta}$
\end{itemize}
}
\end{Definition}

In the decompositions that follow, our trees of models will have the additional property that
$\tp(M_\nu/M_\eta)\perp M_{\eta^-}$
for every $\eta\neq\<\>$ and every $\nu\in Succ_I(\eta)$,
but our early Lemmas do not require this property.

\begin{Lemma} \label{finitetree}
Suppose $\{M_\eta:\eta\in I\}$ is any independent tree of
models indexed by a finite tree $(I,\trianglelefteq)$.  Then the set
$\bigcup_{\eta\in I} M_\eta$ is essentially finite with respect to any strong type
$p$ that is orthogonal to every $M_\eta$.  
\end{Lemma}

\bp  We argue by induction on $|I|$.  For $|I|=1$, this is immediate by Lemma~\ref{basicorth}(1) (taking $A=M_{\<\>}$ and $B=\emptyset$).
So assume $\{M_\eta:\eta\in I\}$ is any independent tree of
models with
$|I|=n+1$ and we have proved the Lemma when $|I|=n$.
Fix any strong type $p$ that is orthogonal to every $M_\eta$.  
Choose any leaf $\eta\in I$ and let $J\subseteq I$ be the subtree with universe
$I\setminus\{\eta\}$.  By the inductive hypothesis, $\bigcup_{\nu\in J} M_\nu$
is essentially finite with respect to $p$,  so  the result follows by
Lemma~\ref{basicorth}(2), taking $A=\bigcup_{\nu\in J} M_\nu$ and
$B=M_\eta$.
\endproof

\begin{Lemma}  \label{anytree}  
Suppose $\{M_\eta:\eta\in I\}$ is any independent tree of
models indexed by any tree $(I,\trianglelefteq)$ and let $N$ be any model 
that contains and is atomic over $\bigcup_{\eta\in I} M_\eta$.
Let $p\in S(N)$ be any regular type that is not eni, but is orthogonal to
every $M_\eta$,  and let $c$ be any realization of $p$.
Then $N\cup\{e\}$ is atomic over
$\bigcup_{\eta\in I} M_\eta$.
\end{Lemma}

\bp  
As notation, for $K\subseteq I$, we let $M_K$ denote $\bigcup_{\nu\in K} M_\nu$.
It suffices to show that $\tp(Dc/M_I)$ is isolated for any finite subset $D\subseteq N$ on which $p$ is
based and stationary.  To see this, fix such a set $D$.  As $D$ is atomic over $M_I$, we can find
a finite set $E\subseteq M_I$ such that $\tp(D/E)$ is isolated and $\tp(D/E)\vdash\tp(D/M_I)$.
Choose a non-empty finite subtree $J\subseteq I$ such that $E\subseteq M_J$ and choose a prime model
$N_J\preceq N$ over $M_J$ that contains $D$.  By Lemmas~\ref{finitetree} and \ref{atomicextension}
we have that $N_J\cup\{c\}$ is atomic over $M_J$.  Choose a formula $\delta(x,h)\in\tp(Dc/M_J)$ 
that isolates the type.  Now, let
$$\F=\{K:J\subseteq K\subseteq I:K\ \hbox{is a subtree and}\ \tp(Dc/M_K) \ \hbox{is isolated by}\ 
\delta(x,h)\}$$
Clearly, $J\in\F$ and by Lemma~\ref{retain} $\F$ is closed under unions of increasing chains.
So choose a maximal element $K^*\in \F$ with respect to inclusion.
To complete the proof of the Lemma, it suffices to prove that $K^*=I$.
If this were not the case, then choose a $\trianglelefteq$-minimal element $\eta\in I\setminus K^*$
and let $K':=K^*\cup\{\eta\}$.
As $J$ was non-empty, $\eta\neq\<\>$ and the independence of the tree yields $\fg {M_{K^*}} {M_{\eta^-}} {M_\eta}$.
But then, by Lemma~\ref{retain}(2), 
$\delta(x,h)$ isolates $\tp(Dc/M_{K'})$, contradicting the maximality of $K^*$.  Thus, $K^*=I$ and the proof is complete.
\endproof

We define a plethora  of decompositions.  

\begin{Definition}  \label{decompdef}
{\em  
Fix a model $M$.  A {\em [regular, eni, eni-active] decomposition inside $M$\/} 
$\d=\<M_\eta,a_\eta:\eta\in I\>$
consists of an independent tree $\{M_\eta:\eta\in I\}$ 
of elementary submodels of $M$
indexed by $(I,\trianglelefteq)$
satisfying the following conditions for each $\eta\in I$:

\begin{enumerate}

\item Each $a_\eta\in M_\eta$ (but $a_{\<\>}$ is meaningless);
\item  The set $C_\eta:=\{a_\nu:\nu\in Succ_I(\eta)\}$ is independent over $M_\eta$;
\item  For each $\nu\in Succ_I(\eta)$ we have:
\begin{enumerate}
\item   $\tp(a_\nu/M_{\nu^-})$ is [regular, eni, eni-active];
\item  If $\eta\neq\<\>$, then $\tp(a_\nu/M_\eta)\perp M_{\eta^-}$; 
\item  $M_\nu$ is atomic over $M_\eta\cup\{a_\nu\}$;
\end{enumerate}
\end{enumerate}

A {\em [regular, eni, eni-active] decomposition of $M$\/} is a
[regular, eni, eni-active] decomposition
inside $M$ with the additional property that for each $\eta\in I$,
the set $C_\eta$ is a {\bf maximal} $M_\eta$-independent set of realizations of [regular, eni, eni-active]
types (that are orthogonal to $M_{\eta^-}$ when $\eta\neq\<\>$).

We say that a decomposition (of any sort) is {\em prime\/} if $M_{\<\>}$ is a prime submodel of $M$
and, for each $\nu\neq\<\>$, $M_{\nu}$ is prime over $M_{\nu^-}\cup\{a_\nu\}$.
}
\end{Definition}

It is important to note that even though eni-NDOP implies eni-active NDOP, it is not the case
that every eni-active decomposition is an eni decomposition.  
As well, note that if $\d=\<M_\eta,a_\eta:\eta\in I\>$ is a decomposition of $M$
(in any of the senses)
and $N\preceq M$ contains $\bigcup_{\eta\in I} M_\eta$, then $\d$ is also a
decomposition of $N$.
The following Lemma requires no assumption beyond $\aleph_0$-stability.

\begin{Lemma}  For any $M$, prime [regular, eni, eni-active]
decompositions of $M$ exist.
\end{Lemma}

\bp  Simply start with an arbitrary prime model $M_{\<\>}\preceq M$, and given a node
$M_\eta$, choose $C_\eta$ to be any maximal $M_\eta$-independent subset of $M$ of realizations of
[regular, eni, eni-active] types (that are orthogonal to $M_{\eta^-}$ when $\eta\neq\<\>$) and, for each
$a_\nu\in C_\eta$, choose $M_\nu\preceq M$ to be prime over $M_\eta\cup\{a_\nu\}$.
Any maximal construction of this sort will produce a prime [regular, eni, eni-active] decomposition of $M$.
\endproof

Of course, without any additional assumptions, such a decomposition may be of limited utility.




\begin{Lemma}[$T$ $\aleph_0$-stable with eni-NDOP] \label{overatomic}
  Let $\d=\<M_\eta,a_\eta:\eta\in I\>$
be any regular decomposition inside $\C$  and let $N$ be atomic over $\bigcup_{\eta\in I} M_\eta$.
If an eni-active regular type $p\not\perp N$, then $p\not\perp M_\eta$ for some $\eta\in I$.
\end{Lemma}

\bp  Recall that eni-NDOP implies eni-active NDOP by
Theorem~\ref{equiv}.  
We first prove the Lemma for all finite index trees $(I,\trianglelefteq)$ by induction on $|I|$.
To begin, if $|I|=1$, then we must have $N=M_{\<\>}$ and there is nothing to prove.
Assume the Lemma holds for all trees of size $n$ and let
  $\d=\<M_\eta,a_\eta:\eta\in I\>$ be a decomposition inside $\C$ indexed by $(I,\trianglelefteq)$ of size
$n+1$.  Let $N$ be atomic over 
$\bigcup_{\eta\in I} M_\eta$ and let $p$ be an eni-active type non-orthogonal to $N$.  
Choose 
a leaf $\eta\in I$ and let $J=I\setminus\{\eta\}$.  
If $(I,\trianglelefteq)$ were a linear order, then again $N=M_\eta$ and there is nothing to prove.
If $(I,\trianglelefteq)$ is not a linear order, then 
choose any $N_J\preceq N$ to be prime over
$\bigcup_{\nu\in J} M_\nu$.
Then, by eni-active NDOP and Lemma~\ref{DOPwitness}(2), either 
$p\not\perp M_\eta$ or $p\not\perp N_J$.  In the first case we are done,
and in the second we finish by the inductive hypothesis since $|J|=n$.
Thus, we have proved the Lemma whenever the indexing tree $I$ is finite.

For the general case, fix a regular decomposition 
  $\d=\<M_\eta,a_\eta:\eta\in I\>$
inside $\C$, let $N$ be atomic over $\bigcup_{\eta\in I} M_\eta$ and choose an eni-active $p\not\perp N$.
By employing Fact~\ref{Fact}(2) and the fact that eni-active types are preserved under non-orthogonality,
we may assume $p\in S(N)$.
Choose a finite $D\subseteq N$ over which $p$ is based
and stationary.  As $D$ is finite and atomic over $\bigcup_{\eta\in I} M_\eta$,
we can find a finite subtree $J\subseteq I$ such that  $D$ is atomic
over $\bigcup_{\eta\in J} M_\eta$.
Fix such a $J$ and choose $M_J\preceq N$ such that $D\subseteq M_J$
and $M_J$ is prime over $\bigcup_{\eta\in J} M_\eta$.
As $D\subseteq M_J$,  $p\not\perp M_J$, so since $J$ is finite, the argument above implies that
there is an $\eta\in J$
such that $p\not\perp M_\eta$.  
\endproof

\begin{Theorem}[$T$ $\aleph_0$-stable with eni-NDOP]  \label{regulartheorem}
Suppose $\<M_\eta,a_\eta:\eta\in I\>$ is a regular decomposition of $M$.
Then $M$ is atomic over $\bigcup_{\eta\in I} M_\eta$.
\end{Theorem}

\bp  Choose $N\preceq M$ to be maximal atomic over 
$\bigcup_{\eta\in I} M_\eta$.  We argue that $N=M$.  If this were not the case,
then choose $e\in M\setminus N$ such that $p:=\tp(e/N)$ is regular.
We obtain a contradiction in three steps.

\medskip
\noindent{{\bf Claim 1:}}  $p\perp M_\eta$ for all $\eta\in I$.
\medskip

\bp  Suppose this were not the case.
Choose $\eta\in I$ $\triangleleft$-minimal
 such that $p\not\perp M_\eta$.  Thus, either $\eta=\<\>$ or
$p\perp M_{\eta^-}$.  
By Lemma~\ref{threemodel}, there is an element $e\in M$ such that $\tp(e/M_\eta)$
is regular and non-orthogonal to $p$ (hence orthogonal to $M_{\eta^-}$ if $\eta\neq\<\>$), 
but $\fg e {M_\eta} {N_\alpha}$.
This element $e$ contradicts the maximality of $C_\eta$ in Definition~\ref{decompdef}.
\endproof

\medskip\noindent{{\bf Claim 2:}}  $p$ is dull.
\medskip

\bp  If $p$ were eni-active, then by Lemma~\ref{overatomic}
we would have $p\not\perp M_\eta$ for some $\eta\in I$, contradicting
Claim 1.
\endproof

As $p$ is dull, it is not eni by Proposition~\ref{eninotdull}.
But this, coupled with Claim~1 implies that $N\cup\{e\}$ is atomic over
$\bigcup_{\eta\in I} M_\eta$, which contradicts the maximality of $N$.
Thus, $N=M$ and we finish.
\endproof

\begin{Theorem}[$T$ $\aleph_0$-stable with eni-NDOP]  \label{eniactivetheorem}  
Suppose $\d=\<M_\eta,a_\eta:\eta\in I\>$ is an eni-active decomposition of a model $M$.
If $N\preceq M$ is atomic over 
$\bigcup_{\eta\in I} M_\eta$, then $N\preceq M$ is a dull pair.
Thus, for every $N'$ satisfying
$N\preceq N'\preceq M$, we have that $N$ is an 
$L_{\infty,\aleph_0}$-elementary substructure of $N'$ and $N'$ is an 
$L_{\infty,\aleph_0}$-substructure of $M$.
\end{Theorem}

\bp  Given $M$ and $\d$, choose any $N\preceq M$ atomic over $\bigcup_{\eta\in I} M_\eta$.
To show that $N\preceq M$ is dull, it suffices to show that there is no $e\in M\setminus N$
such that $\tp(e/N)$ is eni-active.
So, by way of contradiction, assume that there were such an $e$.  Let $p:=\tp(e/N)$.
By Lemma~\ref{overatomic}, we can choose an $\trianglelefteq$-minimal $\eta\in I$
such that $p\not\perp M_\eta$.  By Lemma~\ref{threemodel}, there is $c\in M\setminus M_\eta$
such that $q:=\tp(c/M_\eta)$ is non-orthogonal to $p$ and $\fg c {M_\eta} N$.  
As $q$ is eni-active and orthogonal to $M_{\eta^-}$ (when $\eta\neq\<\>$), the element $c$ contradicts
the maximality of $C_\eta$ in Definition~\ref{decompdef}.  Thus, $N\preceq M$ is a dull pair.
The final sentence follows from Lemma~\ref{dullpairlemma} and Proposition~\ref{DULLchain}.
\endproof

\begin{Corollary}[$T$ $\aleph_0$-stable with eni-NDOP] \label{quotelater}
Suppose $\<M_\eta,a_\eta:\eta\in I\>$ is an eni-active decomposition of both 
$M_1$ and $M_2$.  Then $M_1\bfequiv M_2$.
\end{Corollary}

\bp  Choose any $N_1\preceq M_1$ to be prime over $\bigcup_{\eta\in I} M_\eta$.
By Theorem~\ref{eniactivetheorem}, $N_1\bfequiv M_1$.  By the uniqueness of
prime models, there is $N_2\preceq M_2$ that is both isomorphic to 
$N_1$ and prime over $\bigcup_{\eta\in I} M_\eta$.
By Theorem~\ref{eniactivetheorem} again, $N_2\bfequiv M_2$
and the result follows.
\endproof

The third theorem of this section involves eni decompositions of a model.
Theorem~\ref{enitheorem} is of less interest to us, since when
$M^*$ is uncountable, each of the component submodels $M_\eta$
may be uncountable as well.

\begin{Definition} \label{maxatomic}
 {\em  A decomposition $\<M_\eta,a_\eta:\eta\in I\>$
inside $M$
is {\em maximal atomic} if $M_{\<\>}$ is a maximal atomic substructure
of $M$ and, for each $\nu\neq\<\>$, $M_\nu$ is  maximal atomic
over $M_{\nu^-}\cup\{a_\nu\}$.
}
\end{Definition}

\begin{Theorem}[$T$ $\aleph_0$-stable with eni-NDOP]  \label{enitheorem}
Every model $M$ is atomic over
$\bigcup_{\eta\in I} M_\eta$  for every maximal atomic, 
eni  decomposition
$\<M_\eta:\eta\in I\>$ of $M$.
\end{Theorem}

\bp  Given a maximal, atomic, eni decomposition
$\{M_\eta:\eta\in I\}$ of a model $M$, choose an enumeration $\<\eta_i:i<\alpha\>$ of $I$ such
that $\eta_i\triangleleft\eta_j$ implies $i<j$.  Note that $\eta_0=\<\>$.
Next, define a
continuous, elementary sequence
$\<N_i:i\le\alpha\>$ of elementary substructures of $M$ satisfying:
\begin{itemize}
\item  $N_0=M_{\<\>}$;
\item  $N_\beta=\bigcup_{i<\beta} N_i$ for every non-zero limit ordinal $\beta\le\alpha$; and
\item  $N_{\beta+1}\preceq M$ is maximal atomic over $N_\beta\cup M_{\eta_\beta}$ whenever $\beta<\alpha$.
\end{itemize}
Using Lemma~\ref{retain}, it follows by induction on $\beta\le\alpha$
that  each model $N_\beta$ is atomic over $\bigcup_{i<\beta} M_{\eta_i}$.
Thus, it suffices to prove that $N_\alpha=M$.  Suppose that this were not the case. Choose $e\in M\setminus N_\alpha$
such that $p:=\tp(e/N_\alpha)$ is regular.  Choose $i\le\alpha$ least such that $p\not\perp N_i$.
By superstability, either $i=0$ or $i=\beta+1$ for some $\beta<\alpha$.
We argue by cases, arriving at a contradiction in each case.

\medskip
\noindent{\bf Case 1:}  $p\not\perp M_\eta$ for some $\trianglelefteq$-least $\eta\in I$.
\medskip

\bp  By Lemma~\ref{threemodel}, there is $c\in M\setminus M_\eta$ such that $q:=\tp(c/M_\eta)$ is strongly regular,
non-orthogonal to
$p$, and $\fg c {M_\eta} {N_\alpha}$.  If $q$ were eni, then the element $c$ contradicts the maximality of $C_\eta$
in Definition~\ref{decompdef}.
So assume $q$ is not eni.  There are two subcases:
First, if $\eta=\<\>$, then by Lemma~\ref{atomicextension} (with $A=\emptyset)$ we would have $M_{\<\>}\cup\{c\}$ atomic, contradicting
the maximality of $M_{\<\>}$.  On the other hand, if $\eta\neq\<\>$, then $M_\eta$ would be atomic over $M_\nu\cup\{a_\eta\}$,
where $\nu=\eta^-$.  But then, by Lemma~\ref{basicorth}(1), we would have $M_\nu\cup\{a_\eta\}$ essentially finite with respect to $q$,
hence again by Lemma~\ref{atomicextension} we would have $M_\eta\cup\{c\}$ atomic over $M_\nu\cup\{a_\eta\}$, contradicting the maximality
of $M_\eta$.
\endproof

\medskip\noindent{Case 2:}  $p\perp M_\eta$ for every $\eta\in I$.
\medskip

\bp  In this case, $p$ cannot be dull, because if it were,
then by Lemma~\ref{anytree} $N_\alpha\cup\{e\}$ would be atomic over
$\bigcup_{\eta\in I} M_\eta$.  So assume $p$ is eni-active.
As $N_0=M_{\<\>}$ and $p\not\perp N_i$, the conditions of the Case imply that $i=\beta+1$,
so $N_i$ is atomic over $N_\beta\cup M_\beta$.  Let $\nu=\eta_\beta^-$.  As $N_\beta$ is atomic over
$\bigcup_{\eta_j:j<\beta} M_j$ we have $\fg {N_\beta} {M_\nu} {M_{\eta_\beta}}$.  
Since $p$ is eni-active, then by eni-active NDOP we would have $p\not\perp N_\beta$
or $p\not\perp M_{\eta_\beta}$.  The first possibility contradicts the minimality of $i$, while the second
contradicts the conditions of Case 2.  
\endproof

We close this section with an application of Theorem~\ref{regulartheorem}.
The main point of the proof of  Corollary~\ref{notop} is that
models that are atomic over an independent tree 
of countable models have a large number of partial
automorphisms.

\begin{Corollary}  \label{notop}
If $T$ is $\aleph_0$-stable and eni-NDOP, then $T$ cannot have OTOP.
\end{Corollary}

\bp  By way of contradiction suppose that there were sufficiently
large cardinal $\kappa$ and a model $M^*$ containing a sequence
$\<(b_\alpha,c_\alpha):\alpha<\kappa\>$ and a type $p(x,y,z)$
such that for all $\alpha,\beta<\kappa$,
$$ M^*\ \hbox{realizes}\ p(x,b_\alpha,c_\beta) \quad\hbox{if and only if}
\quad \alpha<\beta$$
For each pair $\alpha<\beta$, fix a realization $a_{\alpha,\beta}$ of
$p(x,b_\alpha,c_\beta)$.
Choose a  prime, regular  decomposition $\<M_\eta,a_\eta:\eta\in I\>$ of $M^*$.
Note that each of the models $M_\eta$ is countable.
By Theorem~\ref{regulartheorem}, $M^*$ is atomic over
$\bigcup_{\eta\in I} M_\eta$, so for each pair $\alpha<\beta$
we can choose a finite $e_{\alpha,\beta}$ from $\bigcup_{\eta\in I} M_\eta$
such that $\tp(a_{\alpha,\beta}/b_\alpha,c_\beta\bigcup_{\eta\in I} M_\eta)$
is isolated by a formula $\theta(x,b_\alpha,c_\beta,e_{\alpha,\beta})$.
We will eventually find a pair $\alpha<\beta$ and $e^*$ from
$\bigcup_{\eta\in I} M_\eta$ such that
$$\tp(b_\beta,c_\alpha,e^*)=\tp(b_\alpha,c_\beta,e_{\alpha,\beta})$$
This immediately leads to a contradiction, as $\theta(x,b_\beta,c_\alpha,e^*)$
would be realized in $M^*$ and any realization of it also realizes
$p(x,b_\beta,c_\alpha)$, contrary to our initial assumptions.

We will obtain these $\alpha<\beta$ and $e^*$ by successively passing from
our sequence to sufficiently long subsequences, each time adding some amount
of homogeneity.
First, for each $\alpha$, choose a finite subtree $J_\alpha\subseteq I$
such that $\tp(b_\alpha c_\alpha/\bigcup_{\eta\in I} M_\eta)$ does not
fork and is as stationary as possible over $J_\alpha$.  By an argument
akin to the $\Delta$-system lemma, by passing to a subsequence we may
assume that there is an $\eta^*\in I$ such that
$J_\alpha\cap J_\beta=\{\nu:\nu\trianglelefteq\eta\}$ for all $\alpha\neq\beta$.
For each $\alpha$, let $M^J_\alpha$ be the countable set
$\bigcup_{\gamma\in J_\alpha} M_\gamma$.
As well,  let $\nu_\alpha$ be the (unique) immediate successor of
$\eta^*$ contained in $J_\alpha$, let
$H_\alpha=\{\gamma\in I:\nu_\alpha\trianglelefteq\gamma\}$,
and let $M_\alpha=\bigcup_{\gamma\in H_\alpha} M_\gamma$.
Note that the sets $H_\alpha$ are pairwise disjoint and the independence of
the tree implies that the sets $\{M_\alpha:\alpha\in\kappa\}$ are independent
over $M_{\eta^*}$.  By trimming further, we may additionally assume that
each of the $J_\alpha$'s are tree isomorphic over $\eta^*$, and that
the sets $M_\alpha$ are isomorphic over $M_{\eta^*}$.

Next, for each $\alpha<\beta$, partition each sequence $e_{\alpha,\beta}$
into three subsequences $r_{\alpha,\beta}\subseteq M_\alpha$,
$s_{\alpha,\beta}\subseteq M_\beta$,
and $t_{\alpha,\beta}$ disjoint from $M_\alpha\cup M_\beta$.

By the Erd\"os-Rado Theorem, we can pass to a subsequence such that
for all $\alpha<\beta<\gamma$ we have:
\begin{itemize}
\item  The partitions coincide, i.e., for each $i$, the $i^{th}$ coordinate
of $e_{\alpha,\beta}\in r_{\alpha,\beta}$ iff the $i^{th}$ coordinate of
$e_{\beta,\gamma}\in r_{\beta,\gamma}$;
\item  $\tp(t_{\alpha,\beta}/M_{\eta^*})$ is constant;
\item  $\tp(r_{\alpha,\beta}/M^J_\alpha)$ is constant; and
\item  $\tp(s_{\alpha,\beta}/M^J_\beta)$ is constant.
\end{itemize}
Additionally, by trimming the sequence still further, we may insist
that for all pairs $\alpha<\beta$, there is
$r^*\in H_\beta$ such that
$\tp(r_{\alpha,\beta}M^J_\alpha/M_{\eta^*})=\tp(r^* M^J_\beta/M_{\eta^*})$
and there is $s^*\in H_\alpha$ such that
$\tp(s_{\alpha,\beta} M^J_\beta/M_{\eta^*})=\tp(s^*M^ J_\alpha/M_{\eta^*})$.

Finally, fix any such $\alpha<\beta$.  By independence, we have
$$\tp(M^J_\alpha,M^J_\beta,r_{\alpha,\beta},s_{\alpha,\beta}, t_{\alpha,\beta})=
\tp(M^J_\beta, M^J_\alpha,r^*,s^*,t_{\alpha,\beta})$$
Let $e^*$ be the sequence formed from $r^*s^*t_{\alpha,\beta}$.
As each of $b_\alpha$ and $c_\beta$ are dominated by $M^J_\alpha$
and $M^J_\beta$, respectively over $M_{\eta^*}$,  it follows
that $\fg {b_\alpha c_\beta} {M^J_\alpha M^J_\beta} {\bigcup_{\eta\in I} M_\eta}$, 
so $\tp(b_\alpha,c_\beta,e_{\alpha,\beta})=\tp(b_\beta,c_\alpha,e^*)$,
completing the proof.
\endproof

\section{Borel completeness of eni-NDOP, eni-deep theories}

Throughout this section, we assume that $T$ is $\aleph_0$-stable with eni-NDOP,
hence prime, eni-active decompositions exist for any model $N$ of $T$.
We begin with a definition, which should be thought of as describing a potential `branch'
of an eni-active decomposition.

\begin{Definition}  \label{CHAIN} 
{\em  An {\em eni-active chain} is a sequence $\<M_i,a_i:i<\alpha\>$, where $2\le\alpha\le\omega$ such that each $a_i\in M_i$
and, for each $i$ such that $(i+1)<\alpha$, $\tp(a_{i+1}/M_i)$ is eni-active,
$\perp M_{i-1}$ (when $i>0$), and $M_{i+1}$ is prime over $M_i\cup\{a_{i+1}\}$.  An eni-active chain is {\em finite} when $\alpha<\omega$.
For $q$ a stationary, regular type, we say a finite chain is {\em $q$-topped\/}
if $q\not\perp M_{\alpha-1}$, but $q\perp M_{\alpha-2}$.  A  finite chain is
{\em ENI-topped} if it is $q$-topped for some ENI type $q$.
}
\end{Definition}

\begin{Definition}  {\em \label{enideepdef}
  An $\aleph_0$-stable, eni-NDOP theory is {\em eni-deep}
if an eni-active  $\omega$-chain exists.
}
\end{Definition}  

Clearly, an $\aleph_0$-stable, eni-NDOP theory is eni-deep if and only if an eni-active $\omega$-chain exists.

\begin{Lemma}  \label{shift}  Given any model $M$ and  regular type $p\in S(M)$, if some stationary, regular type $q$ lies directly over $p$, then there is a $q$-topped finite chain
$\<M_i,a_i:i<\alpha\>$ such that $M_0=M$ and $\tp(a_1/M_0)$ realizes $p$.
\end{Lemma}

\bp  Choose an $\aleph_0$-saturated $N\succeq M$, $a$ realizing $p|N$, and an $\aleph_0$-prime model
$N[a]$ over $N\cup\{a\}$ such that $q\not\perp N[a]$, while
$q\perp N$.  Choose a prime model $M(a)\preceq N[a]$ over $M\cup\{a\}$.
As $q\perp N$, $q\perp M$.  There are now two cases.  First, if $q\not\perp M(a)$,
then  the two-element chain $\<M,M(a)\>$ with $a_1=a$ is as desired.
Second, assume that $q\perp M(a)$.  Choose an eni-active decomposition 
$\<M_\eta:\eta\in I\>$ of $N[a]$
with $M_{\<\>}=M(a)$ such that $M_\eta$ is prime over $M_{\eta^-}\cup\{a_\eta\}$ for every $\eta\in I\setminus\{\<\>\}$.
As $q\not\perp N[a]$ while $q\perp M(a)$, we can choose $\eta\neq\<\>$ minimal such that $q\not\perp M_\eta$.  
As $q\perp M_{\eta^-}$,
$$M\preceq M_{\<\>}\preceq M_{\eta|1}\preceq\dots\preceq M_\eta$$
with $a_{\<\>}=a$ and $a_{\ell+1}=a_{\eta|\ell}$ is a $q$-topped finite chain
as required.
\endproof

Under the assumption of eni-NDOP, this leads to another characterization of the eni-active types.

\begin{Proposition}[$T$ $\aleph_0$-stable, eni-NDOP]  \label{newprop}
A stationary, regular
type $p$ is eni-active if and only if either $p$ is ENI or for every model $M$
such that $p\not\perp M$, there is a finite, ENI-topped chain $\<M_i,a_i:i<\alpha\>$ such that $M_0=M$ and $\tp(a_1/M_0)\not\perp p$.
\end{Proposition}

\bp  Let $\P$ denote the class of types satisfying the alleged characterization.
It follows immediately from Lemma~\ref{notsat} that every type in $\P$ is eni-active.  For the converse, $\P$ visibly contains the ENI types and is closed
under non-orthogonality and automorphisms of the monster model.
Thus, it suffices to show that if $q\in\P$ and $q$ lies directly over $p$,
then $p\in\P$.  To see this, choose any model $M$ such that $p\not\perp M$.
Choose a regular $p'\in S(M)$ non-orthogonal to $p$.  As $q$ lies directly over
$p'$ as well,  use Lemma~\ref{shift} to find a $q$-topped finite chain
$\<M_i,a_i:i<\alpha\>$ with $M_0=M$ and $\tp(a_1/M_0)\not\perp p$.  Now, if $q$ is ENI, then this chain witnesses that
$p\in\P$.  On the other hand, if $q\in\P$ but is not ENI, then there is a finite ENI-topped chain $\<N_j,b_j:j<\beta\>$
with $N_0=M_{\alpha-1}$ and $\tp(b_1/N_0)\not\perp q$.  The concatenation of these two finite chains is an ENI chain starting with $M_0=M$ and $\tp(a_1/M_0)
\not\perp p$.
\endproof

\begin{quotation}  {\bf Until the end of the proof of Theorem~\ref{enideepthm},
fix an $\aleph_0$-stable, eni-NDOP theory that is 
eni-deep as witnessed by a specific eni-active $\omega$-chain 
$\<M_i,a_i:i\in\omega\>$.}
\end{quotation}

Under these hypotheses, we aim to prove Theorem~\ref{enideepthm}.
By employing Proposition~\ref{newprop} for each $i\in\omega$, there is an integer
$k=k(i)>i$ and an ENI-topped finite chain $\CC_{k}=\<N_j^{k},b_j^{k}:j\le k\>$
such that for every $j\le i$, $N^k_j=M_j$ and $b^k_j=a_j$.  
As notation, using Fact~\ref{Fact}(2), choose an ENI $q_k\in S(N^k_k)$ satisfying $q_k\perp N^k_{k-1}$.

\medskip

We will use this configuration of ENI-topped  
chains to code arbitrary subtrees of $\T\subseteq\lambda^{<\omega}$ 
into models $M(\T)$ preserving isomorphism in both directions.  The `reverse direction'
i.e., showing that $M(\T_1)\cong M(\T_2)$ implying 
$(\T_1,\trianglelefteq)\cong (\T_2,\trianglelefteq)$
is quite involved and uses a `black box' in the form of Theorem~6.19 of \cite{ShL}.
We begin by recalling a number of definitions that appear there.
As we are concerned with eni-active decompositions, we take $\P$ to be the
class of eni-active types.  As eni-active types are regular,  
a ${\bf P^r}$-decomposition
in the notation of \cite{ShL} is precisely an eni-active decomposition.

\begin{Definition} {\em  Given a tree $I\subseteq Ord^{<\omega}$, a {\em large subtree}
of $I$ is a non-empty subtree $J\subseteq I$ such that for each $\eta\in J$,
$Succ_I(\eta)\setminus J$ is finite.  We say that two trees $I_1$ and $I_2$ are
{\em almost isomorphic\/} if there exist large subtrees 
$J_1\subseteq I_1$ and $J_2\subseteq I_2$ such that $(J_1,\trianglelefteq)\cong
(J_2,\trianglelefteq)$.

A tree $I$ has {\em infinite branching\/} if, for every $\eta\in I$,
$Succ(\eta)$ is either infinite or empty.
If a tree $I$ has infinite branching,
for any integer $k$, we say a node $\eta\in I$ has {\em uniform depth $k$\/} if, for
every maximal branch of $\{\nu\in I:\eta\trianglelefteq\nu\}$ has length exactly $k$.
A node $\eta$ {\em often has unbounded depth\/} if, for every large subtree $J\subseteq I$
with $\eta\in J$, there is an infinite branch in $J$ containing $\eta$.

Suppose $\eta\in I$ and $E_\eta$ is an equivalence relation on $Succ(\eta)$.
Then $\eta$ is an {\em $(m,n)$-cusp} if there are infinite sets $A_m,A_n,B\subseteq Succ(\eta)$ such that
\begin{enumerate}
\item the set $A_m\cup A_n$ is pairwise $E_\eta$-equivalent;
\item  each $\delta\in A_m$ has uniform depth $m$;  
\item  each $\rho\in A_n$ has uniform depth $n$;
and
\item  each $\gamma\in B$ is often unbounded.
\end{enumerate}
A {\em cusp\/} is an $(m,n)$-cusp for some $m\neq n$.
}
\end{Definition}

\begin{Definition} {\em  Suppose $S\subseteq \P$ and $\d=\<M_\eta,a_\eta:\eta\in I\>$
is a $\P$-decomposition.   We say {\em $\d$ supports $S$\/} if, for every
$q\in S$ there is $\eta(q)\in \max(I)\setminus\{\<\>\}$ such that
$q\not\perp M_{\eta(q)}$, but $q\perp M_{\eta(q)^-}$.  If $\d$ supports $S$,
then we let
Field$(S):=\{\eta(q):q\in S\}$ and $I^S:=\{\nu\triangleleft\eta:\eta\in {\rm Field}(S)\}$.
}
\end{Definition}

\begin{Definition} {\em  Fix a subset $S\subseteq\P$, a model $M$,
 and a function $\Phi:\omega\rightarrow\omega$.
We say that an eni-active decomposition $\d=\<M_\eta,a_\eta:\eta\in I\>$ of $M$
is {\em $\P$-finitely saturated\/} if, for every finite $A\subseteq M$ and $p\in S(A)\cap\P$,
there is $\eta\in I$ such that $\tp(a_\eta/M_{\eta^-})\not\perp p$.

The decomposition $\d$ is {\em $(S,\Phi)$-simple\/} if 
\begin{enumerate}
\item  $\d$ is $\P$-finitely saturated;
\item  $\d$ supports $S$ (hence $I^S$ is defined);
\item  For $\mu\in I$, define $E_\mu$ by $E_\mu(\eta,\nu)\Leftrightarrow \tp(a_\eta/M_\mu)=\tp(a_\nu/M_\mu)$;
\item  for all $\eta,\nu\in I^S$
\begin{enumerate} 
\item  if $\eta^-=\nu^-=\mu$, then $E_\mu(\eta,\nu)$;
\item  $Succ_{I^S}(\eta)$ is empty or infinite (hence $I^S$ has infinite branching);
\item  $\eta$ is either of some finite uniform depth or is a cusp;
\item  if $\eta$ is an $(m,n)$-cusp, then $\Phi(m-n)=\lg(\eta)$.
\end{enumerate}
\end{enumerate}
}
\end{Definition}

Theorem 6.19 from \cite{ShL}, which we take as a black box, states:
\begin{Theorem}  \label{blackbox}
Suppose $S\subseteq\P$, a model $M$, and a function $\Phi:\omega\rightarrow\omega$
are given. If $\d_1$ and $\d_2$ are both $(S,\Phi)$-simple decompositions of $M$, then
the trees $I^S_1$ and $I^S_2$ are almost isomorphic.
\end{Theorem}

With our eye on applying Theorem~\ref{blackbox}, we massage the data we
were given at the top of this section.

Let $\U=\{k\in\omega:k=k(i)$ for some $i\}$.  As $\U$ is infinite, by
passing to an infinite subset, we may additionally assume that
if $n<m$ are from $\U$, then $m>2n$.  It follows from this that
for all pairs $n<m$, $n'<m'$ from $\U$, 
$$m-n=m'-n'\qquad\hbox{if and only if}\qquad m=m'\ \hbox{and}\
n=n'$$
Next, it is routine to partition $\U$ into infinitely many infinite sets $V_i$
for which $k>i$ for every $k\in V_i$.  

Fix an integer $i$.  An `$i$-tree' is a subtree of $\omega^{<\omega}$
with a unique `stem' $\{\<0^j\>:j<i\}$ of length $i$.
As an example, for each $k\in V_i$, let
$$I_i(k):=\{\eta\in\omega^{\le k}:\ \hbox{for all $j<i$, if $\lg(\eta)>j$,
then $\eta(j)=0$\}}$$
If $I$ and $J$ are both $i$-trees (say with disjoint universes) the
{\em free join of $I$ and $J$ over $i$, $I\oplus_iJ$,\/}
is the $i$-tree with universe $(I\cup J)/\sim$, where for each $j<i$,
the (unique) nodes of $I$ and $J$ of length $j$ are identified, and
every other $\sim$-class is a singleton.
To set notation, for $n<m$ from $V_i$, let
$I_i(n,m):=I_i(n)\oplus_i I_i(m)$.
We associate an eni-active decomposition 
$$\d(n,m):=\<N_\eta,b_\eta:\eta\in I_i(n,m)\>$$
satisfying:\begin{itemize}
\item  for $\lg(\eta)<i$, $N_\eta=M_i$ and $b_\eta=a_i$;
\item  if $k(\eta)=n $ when $\eta\in I_i(n)$ and $k(\eta)=m$ when
$\eta\in I_i(m)$, then $N_\eta\cong N_{\lg(\eta)}^{k(\eta)}$ and
$\tp(b_\nu/N_{\nu^-})=\tp(b^{k(\nu)}_{\lg(\nu)}/N^{k(\nu)}_{\lg(\nu^-)})$.
\end{itemize}
In particular, as $\d(n,m)$ is a decomposition,
$\{N_\eta:\eta\in I_i(n,m)\}$ form an independent tree of models.

Still with $i$ fixed, choose disjoint, 4-element sets 
$\{n(\delta^+),m(\delta^+),n(\delta^-),m(\delta^-)\}$ from $V_i$
for each $\delta\in\omega^i$
such that $n(\delta^+)<m(\delta^+)$ and $n(\delta^-)<m(\delta^-)$.

Now, for each $\delta\in\omega^{<\omega}$, let 
$\diff(\delta^+)=m(\delta^+)-n(\delta^+)$ and 
$\diff(\delta^-)=m(\delta^-)-n(\delta^-)$.
It follows from our thinness conditions on $\U$ (and the disjointness of
the sets $V_i$) that the set $D=\{\diff(\delta^+),\diff(\delta^-):\delta\in
\omega^{<\omega}\}$ is without repetition.
Let $\Phi:\omega\rightarrow\omega$ be any function such that for every 
$\delta\in\omega^{<\omega}$,
$$\Phi(\diff(\delta^+))=\Phi(\diff(\delta^-))=\lg(\delta)$$
To ease notation, for each
$\delta\in\omega^{<\omega}$,
let $I(\delta^+)=I_i(n(\delta^+),m(\delta^+))$ and
$\d(\delta^+)=\d(n(\delta^+),m(\delta^+))$, with analogous
definitions for $I(\delta^-)$ and $\d(\delta^-)$.

Next, let $I_0:=(\lambda\times\omega)^{<\omega}$.  We denote
elements of $I_0$ by pairs $(\eta,\delta)$.  Note that $\lg(\eta)=\lg(\delta)$
for all $(\eta,\delta)\in I_0$.  
Let $\d_0$ denote the eni-active decomposition
$\<M_{(\eta,\delta)},a_{(\eta,\delta)}:(\eta,\delta)\in I_0\>$,
where $M_{(\eta,\delta)}\cong M_{\lg(\eta)}$ via a map $f_{(\eta,\delta)}$,
and $f_{(\eta,\delta)}(a_{(\eta,\delta)})=a_{\lg(\eta)}$.

With all of the above as a preamble, we are now ready to code
subtrees of $\lambda^{<\omega}$ into models of our theory.

\begin{Theorem}[$T$ $\aleph_0$-stable, eni-NDOP, eni-deep]  \label{enideepthm}
For any $\lambda\ge\aleph_0$, there is a $\lambda$-Borel embedding $\T\mapsto M(\T)$
of subtrees of $\lambda^{<\omega}$ into models of size $\lambda$ satisfying
$$(\T_1,\trianglelefteq)\cong (\T_2,\trianglelefteq)\qquad\hbox{if and only if}\qquad M(\T_1)\cong M(\T_2)$$
\end{Theorem}

\bp  Fix a cardinal $\lambda\ge\aleph_0$.  We describe the map $\T\mapsto M(\T)$.
Fix a subtree $\T\subseteq\lambda^{<\omega}$.
Begin by letting $\delta_0(\T)$ be the eni-active decomposition
formed by beginning with the decomposition $\d_0$
and simultaneously adjoining a copy of $\d(\delta^+)$ to every
node $(\eta,\delta)\in I_0$ for which $\eta\in \T$, as well as
adjoining a copy of $\d(\delta^-)$ to every node $(\eta,\delta)\in I_0$
for which $\eta\not\in \T$.  
Let $I_0(\T)$ denote the index tree of $\d_0(\T)$.
Let $M_0(\T)$ be prime over $\bigcup\{N_\nu:\nu\in I_0(\T)\}$.
For each $\nu\in\max(I_0(\T))$, let $q_\nu\in S(N_\nu)$ be
the ENI-type conjugate to $q_{\lg(\nu)}\in S(N_{\lg(\nu)}^{\lg(\nu)})$
and let $S=\{q_\nu:\nu\in\max(I_0(\T))\}$.
Because of the independence of the tree and the fact that $M_0(\T)$ is prime
over the tree, each $q_\nu$ has finite dimension in $M_0(\T)$.

Next, we recursively construct an elementary chain $\<M_n(\T):n\in\omega\>$
and a sequence $\<\d_n(\T):n\in\omega\>$ as follows.  We have already defined
$M_0(\T)$ and $\d_0(\T)$, so assume $M_n(\T)$ is defined and
$\d_n(\T)$ is an eni-active decomposition of $M_n(\T)$ extending $\d_0(\T)$.
Let $R_n$ consist of all $p\in S(M_n(\T))\cap\P$ satisfying $p\perp S$.
Let $J_n:=\{a_p:p\in R_n\}$ 
be a $M_n(\T)$-independent set of realizations
of each $p\in R_n$.
For each $p\in R_n$, there is a $\triangleleft$-minimal $\eta(p)\in I_n(\T)$
such that $p\not\perp N_{\eta(p)})$.  
Let $N_p$ be prime over $N_{\eta(p)}\cup \{a_p\}$.
Let $\d_{n+1}(\T)$ be the natural extension of $\d_n(\T)$ formed by 
affixing each $N_p$ as an immediate successor of $N_{\eta(p)}$,
and let $M_{n+1}(\T)$ be prime over the independent tree of models
in $\d_{n+1}(\T)$.  

Finally, let $\d(\T):=\bigcup_{n\in\omega} \d_n(\T)$ and let $M(\T)$ be
prime over $\d(\T)$.  As notation, let $I(\T)$ denote the index tree of $\d(\T)$.

The following facts are easily established:
\begin{enumerate}
\item   A type $p\in S(M(\T))\cap\P$ has finite dimension in $M(\T)$ if and only if
$p\not\perp S$;
\item  $\d(\T)$ is $\P$-finitely saturated;
\item  $\d(\T)$ supports $S$ and $I^S(\T)=I_0(\T)$;
\item $I^S(\T)$ is infinitely branching; and
\item  for $\nu\in I^S(\T)$,
\begin{itemize}
\item   $\nu$ is a cusp if and only if $\nu\in I_0$.  In particular,
if $\nu=(\eta,\delta)$ and $\eta\in\T$, then
$\nu$ is an $(m(\delta^+),n(\delta^+))$-cusp, and $\eta\not\in\T$,
then $\nu$ is an $(m(\delta^-),n(\delta^-))$-cusp;
\item  if $\nu\in I_0(\T)\setminus I_0$, then $\nu$ is of uniform finite depth.
\end{itemize}
\end{enumerate}
In particular, $\d(\T)$ is an $(S,\Phi)$-simple decomposition of $M(\T)$.

\medskip

\par\noindent{{\bf Main Claim:}}  If $M(\T_1)\cong M(\T_2)$, then $(\T_1,\trianglelefteq)\cong
(\T_2,\trianglelefteq)$.

\medskip

\bp  Suppose that $f:M(\T_1)\rightarrow M(\T_2)$ is an isomorphism.
Then the image of $\d(\T_1)$ under $f$ is a decomposition of $M(\T_2)$
with index tree $I(\T_1)$.  
As well, $\d_(\T_2)$ is also a decomposition
of $M(\T_2)$ with index tree $I(\T_2)$.  
If, for $\ell=1,2$, we let $S_\ell$ denote the non-orthogonality classes
of ENI types of finite dimension in $M(\T_\ell)$, then as
isomorphisms preserve types of finite dimension, $f(S_1)=S_2$ setwise.
It follows that both $f(\d_1)$ and $\d_2$ are both $(S_2.\Phi)$-simple
decompositions of $M(\T_2)$.  Thus, by
Theorem~\ref{blackbox}, the trees $I_0(\T_1)$ and $I_0(\T_2)$ are almost
isomorphic.  

Fix large subtrees $J_\ell\subseteq I_0(\T_\ell)$ and a tree isomorphism
$h:J_1\rightarrow J_2$.  
Note that for $\ell=1,2$, a node $\nu\in J_\ell$ has uniform depth $k$ in
$J_\ell$ if and only if $\nu$ has uniform depth $k$ in $I_0(\T_\ell)$.
It follows that $h$ maps cusps to cusps, and more precisely,
$(m,n)$-cusps to $(m,n)$-cusps.
Thus, the restriction $h'$ of $h$ to $J_1\cap(\lambda\times\omega)^{<\omega}$
is a tree isomorphism mapping onto $J_2\cap(\lambda\times\omega)^{<\omega}$
that sends $(m,n)$-cusps to $(m,n)$-cusps.
However, as the pairs $(m,n)$ uniquely identify $\delta\in\omega^{<\omega}$
and even $\delta^+$ and $\delta^-$, it follows that
$h'(\eta,\delta)=(\eta^*,\delta)$ for every $(\eta,\delta)\in \dom(h')$.
As well, if we let 
$$P_\ell:=\{(\eta,\delta)\in J_\ell\cap(\lambda\times\omega)^{<\omega}:
\hbox{$(\eta,\delta)$ is a $\delta^+$-cusp}\}$$
then $h'$ maps $P_1$ onto $P_2$ as well.
Recalling that from our construction, $(\eta,\delta)\in P_\ell$ if and only
if $\eta\in \T_\ell$, we have that for every $(\eta,\delta)\in\dom(h')$
$$\hbox{if $h'(\eta,\delta)=(\eta^*,\delta)$, then $\eta\in \T_1$ if and only
if $\eta^*\in \T_2$.}$$

To finish, we recursively construct maps $h^*:\lambda^{<\omega}\rightarrow
\lambda^{<\omega}$ and $\delta^*:\lambda^{<\omega}\rightarrow
\omega^{<\omega}$ satisfying:
\begin{enumerate}
\item  $(\eta,\delta^*(\eta))\in J_1$;
\item  $h^*(\eta)=\eta^*$ if and only if $h'(\eta,\delta^*(\eta))=(\eta^*,
\delta^*(\eta))$;
\item for all $\eta$ and all $\alpha,\alpha'\in\lambda$,
$\delta^*(\eta\conc\<\alpha\>)=\delta^*(\eta\conc\<\alpha'\>)$; and
\item  for all $\eta\in \lambda^{<\omega}$, $\alpha,\beta\in\lambda$,
$(\eta\conc\<\alpha\>,\delta^*(\eta\conc\<\alpha\>))\in J_1$ and \hfill \break
$(h^*(\eta)\conc\<\beta\>,\delta^*(h^*(\eta)\conc\<\beta\>))\in J_2$.
\end{enumerate}
To accomplish this, first let $\delta^*(\<\>)=\<\>$.
Given that $\delta^*(\eta)$ is defined, the definition of $h^*(\eta)$ is
given by Clause~(2).  As $(\eta,\delta^*(\eta)\in J_1$
and since $J_\ell$ are large subtrees of $I_0(\T_\ell)$,
it follows that there is $\delta'\in Succ(\delta^*(\eta))$
such that Clauses (3) and (4) hold for all $\alpha,\beta\in\lambda$.
Define $\delta^*(\eta\conc\<\alpha\>)=\delta'$ for every $\alpha$
and define $h^*(\eta\conc\<\alpha\>)$ according to Clause~(2).

It is easily checked that 
$h^*:\lambda^{<\omega}\rightarrow\lambda^{<\omega}$ is a tree
isomorphism.  Additionally, as $h'$ mapped $P_1$ onto $P_2$,
it follows that the restriction of $h^*$ to $\T_1$ is a tree isomorphism
between $(\T_1,\trianglelefteq)$ and $(\T_2,\trianglelefteq)$.
\endproof

\begin{Corollary}  \label{enideepcor}
If $T$ is $\aleph_0$-stable with eni-NDOP and is eni-deep, then 
$T$ is Borel complete.  Moreover, for every infinite cardinal $\lambda$,
$T$ is $\lambda$-Borel complete for $\bfequiv$.
\end{Corollary}

\bp  If $T$ has eni-DOP, then this is literally Corollary~\ref{eniDOPcor}.
If $T$ has eni-NDOP, then the proof is exactly like the proof of 
Corollary~\ref{eniDOPcor}, using Theorem~\ref{enideepthm} in place of 
Theorem~\ref{eniDOPthm}.
\endproof

\section{Main gap for  models of $\aleph_0$-stable theories modulo $L_{\infty,\aleph_0}$-equivalence}
\label{last}

In this brief section, we  combine our
previous results to exhibit a dichotomy among $\aleph_0$-stable theories.

\begin{Definition}  {\em  For $T$ any theory and $\lambda$ an infinite
cardinal, let $Mod_\lambda(T)$ denote the set of models of $T$ with universe
$\lambda$.
\begin{itemize}
\item  For $T$ any theory and $\lambda$ any cardinal, 
$I_{\infty,\aleph_0}(T,\lambda)$ denotes the maximum cardinality of
any pairwise non-$\bfequiv$ collection from $Mod_\lambda(T)$.
\item  For any $M\models T$ of size $\lambda$,
the {\em Scott height of $M$}, $SH(M)$ is the least ordinal $\alpha<\lambda^+$
such that for any model $N$, $N\equiv_\alpha M$ implies 
$N\equiv_{\alpha+1} M$.
\end{itemize}
}
\end{Definition}

\begin{Theorem}  \label{maingap}  
The following conditions are equivalent for any $\aleph_0$-stable theory $T$:
\begin{enumerate}
\item  For all infinite cardinals $\lambda$,
$I_{\infty,\aleph_0}(T,\lambda)=2^\lambda$;
\item  For all infinite cardinals $\lambda$, 
$\sup\{SH(M):M\in Mod_\lambda(T)\}=\lambda^+$;
\item  $T$ either has eni-DOP or is eni-deep.
\end{enumerate}
\end{Theorem}

\bp  The equivalence of $(1)\Leftrightarrow(2)$ is the content of
\cite{Sh11}.  

$(3)\Rightarrow(1):$ Fix any infinite cardinal $\lambda$.
If $T$ has either of these properties,
then  by Corollary~\ref{eniDOPcor} or Corollary~\ref{enideepcor},
$T$ is $\lambda$-Borel complete.  However, it is well known (see e.g.,
\cite{Sh220})
that there is a family of $2^\lambda$ pairwise non-$\bfequiv$ directed
graphs with universe $\lambda$.  It follows immediately that
$I_{\infty,\aleph_0}(T,\lambda)=2^\lambda$ in either case.

$(1)\Rightarrow(3):$  Assume that $T$ is $\aleph_0$-stable, with
eni-NDOP and eni-shallow (i.e., not eni-deep).  Then, 
by Corollary~\ref{quotelater}, models of $T$ are determined by
up to $\bfequiv$-equivalence by their prime, eni-active decompositions.
Thus, it suffices to count the number of prime, eni-active decompositions
up to isomorphism.\footnote{We say that two eni-active decompositions 
$\d_1=\<M_\eta^1,a^1_\eta:\eta\in I_1\>$ and
$\d_2=\<M_\eta^2,a^2_\eta:\eta\in I_2\>$ are isomorphic if there
is a tree isomorphism $f:(I_1,\trianglelefteq)\cong (I_2,\trianglelefteq)$
and an elementary bijection 
$f^*:\bigcup_{\eta\in I_1} M^1_\eta\rightarrow\bigcup_{\eta\in I_2} M^2_\eta$
such that, for each $\eta\in\ I_1$, $f^*|_{M^1_\eta}$ maps $M^1_\eta$ isomorphically onto $M^2_{f(\eta)}$.}

To obtain this count, first note that if $T$ is eni-shallow, then 
as in Theorem~X~4.4  of \cite{Shc} (which builds on VII, Section~5
of \cite{Shc}), the depth of any index tree of
an eni-active decomposition is an ordinal $\beta<\omega_1$.  
In any prime decomposition, each of the models $M_\eta$ is countable,
hence there are at most $2^{\aleph_0}$ isomorphism types.
So, as a weak upper bound, if $\lambda=\aleph_\alpha$, then the number
of prime, eni-active decompositions of depth $\beta$ 
of a model of size $\lambda$ is bounded
by $\beth_{(|\alpha|+|\beta|)^+}$. [Similar counting arguments
appear in Theorem~X~4.7 of \cite{Shc}.] From this, we conclude that
for some cardinals $\lambda$, 
$I_{\infty,\aleph_0}(T,\lambda)<2^\lambda$.
\endproof

\appendix

\section{Appendix: Packing problems for bipartite graphs}

A bipartite graph $A$ consists of a set of vertices, which are partitioned into
 two sets
$L(A)$ and $R(A)$, together with a binary,  irreflexive edge
relation $E(A)\subseteq L(A)\times R(A)$.
We say that $A$ is {\em complete bipartite}
if the set of edges $E(A)=L(A)\times R(A)$.
We call $A$ {\em balanced} if $||L(A)|-|R(A)||\le 1$.

Define a function $e^*:\omega\rightarrow\omega$ by
$e^*(2b)=b^2$ and $e^*(2b+1)=b(b+1)$ for all $b\in\omega$.
A classical packing problem asserts:

\begin{Fact}
\label{classical} A  bipartite graph $A$ with at most $c\ge 2$ vertices
has at most $e^*(c)$ edges, with equality holding if and only if
$|A|=c$ and $A$ is complete and balanced.
\end{Fact}

For a bipartite graph $A$, let $v(A)$, $e(A)$, and $CC(A)$
denote the number of vertices, edges, and connected components of $A$, respectively.  Recall
that in the discussion prior to the statement of Proposition~\ref{biggy},  we defined
an $m_1\times m_2$ bipartite graph $A$ to be almost $\ell$-complete if
$|m_i-\ell|\le 0.01\ell$ for each $i=1,2$ and each vertex has valence
at least $0.9\ell$.

\begin{Fact}  \label{almostell}  Suppose that $N$ is a given integer
and  $\ell>>N$ (explicit bounds on $\ell$ in terms of $N$ can be
 found from the proof).
 If $A$ is any bipartite graph with $v(A)\le 2\ell+N$
and $e(A)\ge\ell^2-N$, then $A$ is almost $\ell$-complete.
\end{Fact}

\bp
The new statistic we investigate in this Appendix  is $k(A)$, which we define to be $v(A)-CC(A)$.
Two special cases are that $v(A)=k(A)+1$ for any connected bipartite graph $A$,
and that any null bipartite graph $B$ has $k(B)=0$.

We wish to find an analogue of Fact~\ref{almostell} in which the upper
bound on $v(A)$ is replaced by an upper bound on $k(A)$.

If $A$ and $B$ are each bipartite graphs with disjoint sets of vertices, then
$A\coprod B$ denotes their {\em disjoint union}.  It is the bipartite graph $C$  whose
vertices are the union of the vertices of $A$ and $B$, and
$E(C)=E(A)\cup E(B)$.

Note that all of our statistics are additive with respect to disjoint unions.  For example,
for $x\in\{n,e,CC,k\}$, $x(A\coprod B)=x(A)+ x(B)$.  Thus, if $A$ is any bipartite graph
and $B$ is null, then $k(A\coprod B)=k(A)$.
The proof of the following Lemma is routine.

\begin{Lemma}   \label{1}
Suppose $A$ and $B$ are disjoint, and are each complete, balanced, bipartite graphs
with $k(A)\ge k(B)\ge 1$.  Let $A^+$  and $B^-$ be disjoint, complete, balanced bipartite graphs
with $k(A^+)=k(A)+1$ and $k(B^-)=k(B)-1$.   Then $k(A^+\coprod B^-)=k(A\coprod B)$
and $e(A^+\coprod B^-)\ge e(A\coprod B)$.
\end{Lemma}

One Corollary follows immediately by combining Fact~\ref{classical} with
Lemma~\ref{1}.

\begin{Corollary}  For all positive integers $a$ and all bipartite graphs $A$ with $k(A)\le a$,
$e(A)\le e^*(a+1)$, with equality holding if and only if $A=B\coprod C$,
with $B$ complete and balanced,
and $C$ null (and may be empty).
\end{Corollary}

Next, given a pair of integers $c,d$, let $f(c,d)$ be the least integer such that $e(A)\le f(c,d)$
for all bipartite graphs of the form  $A=B\coprod C$, where $k(B)\le c$ and $k(C)\le d$.

\begin{Lemma}  \label{f}
\begin{enumerate}
\item  For all $c,d\in\omega$, $f(c,d)=e^*(c+1)+e^*(d+1)$; and
\item  If $1\le d\le c$, then $f(c+1,d-1)\ge f(c,d)$.
\end{enumerate}
\end{Lemma}

\bp  The first statement follows by applying Fact~\ref{classical} to each
of $B$ and $C$, while the second follows from Lemma~\ref{1}.
\endproof

\begin{Proposition} \label{bound}
If $\ell>W^2/4$ and $A$ is a bipartite graph satisfying $k(A)\le 2\ell+W$ and $e(A)\ge \ell^2$,
then $A$ contains a connected subgraph $B\subseteq A$ with  at least $\ell^2-W^2/4$ edges and at most $2\ell+W$ vertices.
\end{Proposition}

\bp  Let $\Phi$ be the set of all $A$ such that $k(A)\le 2\ell +W$ and $A$ does not have any connected component
$B$ with $k(B)\ge 2\ell-1$.  Among all such $A$,
choose  $A^*\in\Phi$ so as to maximize the number of edges $e(A^*)$.
By Lemma~\ref{f},
$$e(A^*)\le f(2\ell-1,W+1)=e^*(2\ell-1)+e^*(W+1)\le\ell(\ell-1) + W^2/4<\ell^2$$
Thus, our given graph $A\not\in\Phi$,  so $A$ has a connected component $B$
with $k(B)\ge 2\ell-1$.  Now, if we decompose $A$ as $A=B\coprod C$,
then $k(C)\le W+1$, so by Lemma~\ref{1},
$e(C)\le e^*(W+1)\le W^2/4$.

Since $\ell^2\le e(A)=e(B)+e(C)$, this implies $e(B)\ge \ell^2-W^2/4$.
But, since $B$ is connected, $v(B)=k(B)+1$ and $k(B)\le k(A)\le 2\ell+W$,
so $B$ has at most $2\ell+W$ vertices.
\endproof

The following Corollary follows immediately by combining Proposition~\ref{bound}
with Fact~\ref{almostell}.

\begin{Corollary}  \label{finitecomb}
If $\ell>>W$ and $A$ is a bipartite graph satisfying
$k(A)\le 2\ell+W$ and $e(A)\ge\ell^2$, then $A$ has an almost
$\ell$-complete subgraph.
\end{Corollary}


\begin{thebibliography}{99}



\bibitem{BL} E.\ Bouscaren and D.\ Lascar,
Countable models of nonmultidimensional $\aleph_{0}$-stable theories,
{\it Journal of  Symbolic Logic} {\bf 48}  (1983),  no. 1, 197--205.

\bibitem{FS} H.\ Friedman and L.\ Stanley,
A Borel reducibility theory for classes of countable structures,
{\it   Journal of  Symbolic Logic} {\bf 54}  (1989),  no. 3, 894--914.


\bibitem{K} M.\ Koerwien,
Comparing Borel reducibility and depth of an $\omega$-stable theory,
{\it  Notre Dame J.\ of Formal Logic\/} {\bf 50} (2009) 365-380.

\bibitem{L} M.C.\ Laskowski, An old friend revisited:
Countable models of $\omega$-stable theories,
{\it Notre Dame J.\ of  Formal Logic} {\bf 48} (2007) 133-141.

\bibitem{Makkai} M. Makkai, A survey of basic stability theory
with particular emphasis on orthogonality and regular types,
{\it Israel Journal of Mathematics\/}
{\bf 49} (1984), 181-238.

\bibitem{Shc} S. Shelah, {\it Classification Theory\/}, (revised edition)
North Holland, Amsterdam, 1990.

\bibitem{Sh11} S. Shelah, On the number of non-almost isomorphic models of $T$ in a power,
{\it Pacific Journal of  Math} {\bf 36} (1971) 811-818.

\bibitem{Sh220} S. Shelah, Existence of many $L_{\infty,\lambda}$-equivalent, nonisomorphic models of 
$T$ of power $\lambda$, 
{\it Annals Pure and Applied Logic} {\bf 34} (1987) 291-310. 

\bibitem{Sh401} S.\ Shelah,
Characterizing an $\aleph_\epsilon$-saturated model of superstable
NDOP theories by its $\Bbb L_{\infty,\aleph_\epsilon}$-theory,
{\it  Israel Journal of  Math} {\bf  140} (2004) 61-111.

\bibitem{ShL} S.\ Shelah and M.C.\ Laskowski,
{\bf P}-NDOP and {\bf P}-decompositions of $\aleph_\epsilon$-saturated models
of superstable theories, preprint.

\bibitem{SHM} S.\ Shelah, L.A.\ Harrington, and M.\ Makkai,
A proof of Vaught's conjecture for $\omega$-stable theories,
{\it Israel Journal of  Math} {\bf 49}(1984),  no. 1-3, 259--280.

\end{thebibliography}
\end{document}